\documentclass[12pt]{amsart}

\usepackage{amssymb,latexsym}
\usepackage{enumitem}
\usepackage{comment}

\makeatletter

\@namedef{subjclassname@2010}{

  \textup{2010} Mathematics Subject Classification}
\usepackage{amssymb, amsmath,amsthm}
\newtheorem{thm}{Theorem}[section]
\newtheorem{prop}[thm]{Proposition}
\newtheorem{conj}[thm]{Conjecture}

\newtheorem{cor}[thm]{Corollary}
\newtheorem{lem}[thm]{Lemma}
\theoremstyle{definition}
\newtheorem{rem}[thm]{Remark}
\newtheorem{def1}[thm]{Definition}

\newcommand{\ra}{\rightarrow}
\newcommand{\bk}{\backslash}
\newcommand{\mc}{\mathcal}
\newcommand{\mf}{\mathfrak}
\newcommand{\mb}{\mathbb}
\newcommand{\sg}{\sigma}
\newcommand{\eps}{\epsilon}

\newcommand{\llf}{\left\lfloor}

\newcommand{\e}{\varepsilon}
\newcommand{\rrf}{\right\rfloor}

\newcommand{\mbf}{\boldsymbol}
\renewcommand{\bar}{\overline}
\newcommand{\asum}{\sideset{}{^{\ast}}\sum}

\newcommand{\proofpart}[2]{%
  \par
  \addvspace{\medskipamount}%
  \noindent\emph{Part #1: #2}\par\nobreak
  \addvspace{\smallskipamount}%
  \@afterheading
}

\begin{document}
\title[Inequalities with Coefficients of Cusp Forms]{Monotone Chains of Fourier Coefficients of Hecke Cusp Forms}

\author{Oleksiy Klurman}
\address{Max Planck Institute for Mathematics, Bonn, Germany and 
School of Mathematics, University of Bristol, Bristol, UK}
\email{lklurman@gmail.com}

\author{Alexander P. Mangerel}
\address{Centre de Recherches Math\'{e}matiques, Universit\'{e} de Montr\'{e}al, Montr\'{e}al, Qu\'{e}bec}
\email{smangerel@gmail.com}

\begin{abstract}
We prove general equidistribution statements (both conditional and unconditional) relating to the Fourier coefficients of arithmetically normalized holomorphic Hecke cusp forms $f_1,\ldots,f_k$ without complex multiplication, of equal weight, (possibly different) squarefree level and trivial nebentypus.\\
As a first application, we show that  for the Ramanujan $\tau$ function and any admissible $k$-tuple of distinct non-negative integers $a_1,\ldots,a_k$ the set
$$
\{n \in \mb{N} : |\tau(n+a_1)| < \cdots < |\tau(n+a_k)|\}
$$
has positive natural density. This result improves upon recent work of Bilu, Deshouillers, Gun and Luca [{\it Compos. Math.} (2018), no. 11, 2441-2461]. Secondly, we make progress towards understanding the signed version by showing that
$$
\{n \in \mb{N} : \tau(n+a_1) < \tau(n+a_2) < \tau(n+a_3)\}
$$
has positive relative upper density at least $1/6$ for any admissible triple of distinct non-negative integers $(a_1,a_2,a_3).$ More generally, for such chains of inequalities of length $k > 3$ we show that under the assumption of Elliott's conjecture on correlations of multiplicative functions, the relative natural density of this set is $1/k!.$ Previously results of such type were known for $k\le 2$ as consequences of works by Serre and by Matom\"{a}ki and  Radziwi\l\l.\\
Our results rely crucially on several key ingredients: i) a multivariate Erd\H{o}s-Kac type theorem for the function $n \mapsto \log|\tau(n)|$, conditioned on $n$ belonging to the set of non-vanishing of $\tau$, generalizing work of Luca, Radziwi\l\l \ and Shparlinski; ii) the recent breakthrough of Newton and Thorne on the functoriality of symmetric power $L$-functions for $\text{GL}(n)$ for all $n \geq 2$ and its application to quantitative forms of the Sato-Tate conjecture; and iii) the work of Tao and Ter\"{a}v\"{a}inen on the logarithmic Elliott conjecture.
\end{abstract}
 
\maketitle

\section{Introduction}

Let $f: \mb{N} \ra \mb{R}$ be a multiplicative function. 
It is a well-known result of Erd\H{o}s \cite{Erd} that unless $f(n) = n^{\alpha}$ for some real $\alpha \neq 0$ and all $n \in \mb{N}$ (or $f(n) = 0$ for all $n \geq 2$) then $\{f(n)\}_n$ is not monotone. Avoiding these monotone cases, it is generally a very difficult problem to determine how often $f$ is increasing or decreasing along a tuple of consecutive integers. In this direction, we recall a conjecture of  of Sarkozy \cite[Conjectures 53 and 54]{Sar}.
\begin{conj}[Sark\H{o}zy] \label{Sark}
If $f:\mathbb{N}\to\mathbb{N}$ is a multiplicative function such that for any $k\ge 0$ there are at least two distinct primes $p$ for which $f(p)\ne p^k,$ then both of the inequalities 
$$f(n)>\max\{f(n-1),f(n+1)\}$$
and
$$f(n)<\max\{f(n-1),f(n+1)\}$$
hold for infinitely many $n\ge 2.$
\end{conj}
Conjecture~\ref{Sark} is a manifestation of a more general phenomenon, namely that given a typical multiplicative function $f: \mb{N} \ra \mb{R}$ and any admissible set of distinct integers $a_1,a_2\dots, a_k$ the set 
\begin{equation}\label{chaink}
\{n \in \mb{N}, f(n+a_i)\ne 0, 1 \leq i \leq k: f(n+a_1) \leq f(n+a_2) \leq \cdots \leq f(n+a_{k})\}
\end{equation}
should be infinite and moreover have positive density.
We highlight here a very important special case of~\eqref{chaink} that of $f(n)=\text{sign}(g(n))$ for a given multiplicative $g:\mathbb{N}\to\mb{R}.$ In this case, our problem reduces to another widely studied question, that of counting tuples of values of $f$ with prescribed sign patterns. For generic, unbounded multiplicative functions $f$, the event $f(n+a_i) = f(n+a_j)$ is rare, and so for the purposes of this paper we will count sets like in \eqref{chaink}, where the inequalities are to be replaced by strict ones.

One of the objectives of the present paper is to study~\eqref{chaink} for multiplicative functions determined by the Fourier coefficients of non-CM primitive Hecke cusp forms. To this end, we
let $f(z) := \sum_{n \geq 1} b_f(n)e(nz)$ be a primitive non-CM Hecke eigencusp form\footnote{Throughout this paper, whenever we refer to a primitive cusp form, we mean a primitive Hecke eigenform that is normalized so that its first Fourier coefficient is 1.} weight $m$, square-free level $N$ and trivial nebentypus, normalized so that $b_f(1) = 1$, and let $\lambda_f(n) := b_f(n)n^{-(m-1)/2}.$ It is well-known that $\{b_f(n)\}_n$ is a sequence of totally real algebraic numbers, and so $n \mapsto b_f(n)$ is a real-valued multiplicative function. To exemplify some of our main results, we recall that if $f = \Delta$, the unique cusp form of weight $m=12$ for the full modular group (i.e., $N=1$), then $b_f(n):=\tau(n),$ is the well-known Ramanujan $\tau$ function, which is integer-valued. Our guiding conjecture in this direction is the following folklore statement. 
\begin{conj} \label{conj:serre}
Let $k \geq 1$ and let $a_1,\ldots,a_k$ be distinct non-negative integers for which $\tau(a_j) \neq 0$. Then the set
\begin{equation}\label{tauchain}
\{n \in \mb{N}, \tau(n+a_i)\ne 0, 1 \leq i\leq k : \tau(n+a_1) <\tau(n+a_2) < \cdots < \tau(n+a_{k})\}
\end{equation}
is infinite.
\end{conj}
A famous conjecture of Lehmer states that $\tau(n)\ne 0$ for all $n\ge 2$. Under this assumption, it is natural to upgrade the statement of~\eqref{conj:serre} and conjecture that the natural density (the definition of which is recalled below) of the corresponding set exists and is in fact equal to $\frac{1}{k!}.$  This would reflect the order statistics that one might expect from a random sequence of real numbers, and a number of papers have borne witness to the pseudorandom behaviour of the sequence $\{\tau(n)\}_n$ and of the sequences of Fourier coefficients of other Hecke cusp forms, as we discuss in Section \ref{sssec:signs} below.
\subsubsection{\bf Unsigned orderings of primitive cusp forms.} 
One of the main obstacles in effectively studying~\eqref{tauchain} is the oscillating sign of $\tau.$ To avoid this difficulty, one can therefore consider the corresponding question for $|\tau(n)|.$ As an application of the main theorem in the recent beautiful paper \cite{BDGL} (see Theorem 1.1 there), the authors show that if $a_1,\dots,a_{k}$ are distinct non-negative integers such that $\tau(a_j) \neq 0$ then 
$$
\{n \in \mb{N} : 0 < |\tau(n+a_1)| < |\tau(n+a_2)| < \cdots < |\tau(n+a_{k})|\}
$$
contains infinitely many integers, and moreover the count of such $n \leq x$ is $\gg_k x/(\log x)^k$. Our first theorem provides an improvement of this result.
Recall that the \emph{natural density} of a set $S \subseteq \mb{N}$, is defined by
$$
d(S) := \lim_{Y \ra \infty} Y^{-1}|\{n \leq Y : n \in S\}|,
$$
whenever this limit exists. 
\begin{cor} \label{cor:ramTau}
Let $k \geq 1$ and let $a_1,\ldots,a_k$ be distinct non-negative integers for which $\tau(a_j) \neq 0$. Then the set
$$
\{n \in \mb{N}: 0 < |\tau(n+a_1)| < |\tau(n+a_2)| < \cdots < |\tau(n+a_{k})|\}
$$
has positive density. Moreover, if Lehmer's conjecture is true (i.e., $\tau(n) \neq 0$ for all $n \in \mb{N}$) then the density is $1/k!$.
\end{cor}
\begin{rem} We remark that in~\cite{BDGL}, a suitable set of $n\le x$ is constructed such that $\omega(\prod_{i\le k}(n+a_i))\ll_k 1,$ and consequently their number is necessarily sparse, specifically $\ll_k x/(\log x)^{k-o(1)}.$ In our case, the integers $n\le x$ in the conclusion of Corollary~\ref{cor:ramTau} satisfy $\omega(n+a_i)\sim \log\log x$ for all $1\le i\le k,$ which naturally provides a density increase.
\end{rem}
The above result follows from a more general statement, applying to the Fourier coefficients of a collection of possibly different primitive non-CM cusp forms having the same weight, possibly different squarefree levels and all of trivial nebentypus. To state it more precisely, we require the following definition.
\begin{def1}
Let $f_1,\ldots,f_k$ be primitive non-CM cusp forms with respective sequences of Fourier coefficients $\{b_j(n)\}_n$. A $k$-tuple $\mbf{a} = (a_1,\ldots,a_k)$ is said to be \emph{$\mbf{f}$-admissible} if, for every prime $p$ such that $p|\prod_{1 \leq i < j \leq k} (a_j-a_i)$ and $b_1(p)\cdots b_k(p) = 0$, the set 
$$
\{c \pmod{p} : \prod_{1 \leq j \leq k} (c+a_j) \not \equiv 0 \pmod{p}\}
$$ 
is non-empty.
\end{def1}
\begin{rem}
It is known that $\tau(n)\ne 0$ for all $n\le 10^{22}$ and therefore all the finite subsets of $[1,10^{22}]$ are $\mbf{\tau}$-admissible.
\end{rem}
\begin{thm} \label{thm:genCusp}
Let $k\geq 1$. Let $f_1,\ldots,f_k$ be arithmetically normalized primitive Hecke cusp forms of the same weight, and each of squarefree level and trivial nebentypus, and let $a_1,\ldots,a_k$ be distinct non-negative integers such that $\mbf{a}$ is an $\mbf{f}$-admissible $k$-tuple. Write $f_j(z) = \sum_{n \geq 1} b_j(n)e(nz)$, and assume that for each $1 \leq j \leq k$ the following assumption holds:
\begin{enumerate}[label = (\Alph*)]
\item if $b_j(p) \neq 0$ then $b_j(p^{\nu}) \neq 0$ for all primes $p$ and integers $\nu \geq 1$, and moreover
$|b_j(p^{\nu})| \geq p^{-C\nu}$ for some constant $C = C(f_j) > 0$. 
\end{enumerate}
Then
$$
\{n \in \mb{N}: 0 <|b_1(n+a_1)| < \cdots < |b_k(n+a_k)|\}
$$
has positive density. In fact, if we set 
$$
\mc{N} := \{n \in \mb{N} : \prod_{1 \leq j \leq k} b_j(n) \neq 0\},
$$
we then have
$$
\lim_{x \ra \infty} \frac{|\{n \leq x : n+a_1,\ldots,n+a_k \in \mc{N}, 0 < |b_1(n+a_1)| < \cdots < |b_k(n+a_k)|\}|}{|\{n \leq x : n+a_1,\ldots,n+a_k \in \mc{N} \}|} = \frac{1}{k!}.
$$
\end{thm}
This improves upon the main Theorem 1.2 in \cite{BDGL}, where, as mentioned above, it is shown that for a number $\gg_k x/\log^k x$ of positive integers $n \leq x$ the inequalities $|b_j(n+a_j)| < |b_{j+1}(n+a_{j+1})|$ hold for all $1 \leq j \leq k-1$.
\begin{rem} \label{rem:NjPD}
It can be shown that $\mc{N}$ is a positive density set, and in fact so is 
$$
\{n \in \mb{N} : n+a_j \in \mc{N} \text{ for all } 1\leq j \leq k\}
$$ 
whenever $\mbf{a}$ is $\mbf{f}$-admissible (see Proposition \ref{prop:Sieve} below). Thus, the limit in the statement is well-defined. This can be seen as a generalization of a result of Serre \cite{Ser} who proved the corresponding statement for $k=1.$
\end{rem}
\begin{rem}
It is known that any non-CM cusp form $f$ satisfies the property $b_j(p) \neq 0 \Rightarrow b_f(p^{\nu}) \neq 0$ for all $\nu$ for all but finitely many $p$. As noted on p. 175 of \cite{LRS}, the lower bound $b_f(p^{\nu}) \geq p^{-C\nu}$ for some $C> 0$ holds for all holomorphic Hecke cusp forms, which includes $f = \Delta$. This is why no assumptions are required for the Ramanujan $\tau$ function.
\end{rem}
\subsubsection{\bf Signed orderings of cusp forms.} \label{sssec:signs}
Roughly speaking, our treatment of the monotone patterns \eqref{chaink} amounts to understanding both the joint distribution of the magnitude vector $(|\tau(n+a_1)|,|\tau(n+a_2)|, \ldots,|\tau(n+a_k)|)$ together with information about the patterns of the vector $(\text{sign}(\tau(n+a_1)),\text{sign}(\tau(n+a_2)), \ldots,\text{sign}(\tau(n+a_k)))$, which captures the oscillations of signs. 

In~\cite{KLSW}, Kowalski, Lau, Soundararajan and Wu investigated the problem of bounding the first sign change of $\lambda_f$ for holomorphic $f$, and showed that the sequence of signs $\{\text{sign}(\lambda_f(p))\}_p$ determines the underlying form $f$. The frequency of sign changes and related questions have been extensively studied in recent years (see e.g., \cite{K}, \cite{KLSW} and \cite{GS}) culminating in the work of Matom\"{a}ki and Radziwi\l\l~\cite{MR}. They showed that the sequence $\{\lambda_f(n)\}_n$ has a positive proportion of sign changes for both holomorphic cusp forms and Hecke-Maass forms for the full modular group. In the former case this result later became a direct consequence of their landmark work~\cite{MR3488742}.

In our modest attempt towards resolving Conjecture~\ref{conj:serre},  we state both conditional and unconditional results, each of the flavour that sets like that in \eqref{chaink} have either positive natural or upper density. We recall that the \emph{upper (natural) density} of a set $S \subset \mb{N}$ is defined as
$$
\bar{d}(S) := \limsup_{Y \ra \infty} Y^{-1} |\{n \leq Y : n \in S\}|.
$$
Our unconditional result is as follows.
\begin{thm} \label{thm:uncondSigns}
Let $k \in \{2,3\}$ and let $a_1,\ldots,a_k$ be distinct non-negative integers. Let $f_1,\ldots,f_k$ be arithmetically normalized primitive non-CM Hecke cusp forms of equal weight, squarefree level and trivial nebentypus, satisfying the conditions in Theorem \ref{thm:genCusp}, with respective sets of Fourier coefficients $\{b_j(n)\}_n$ for $j = 1,\ldots,k$. Suppose also that $\mbf{a}$ is $\mbf{f}$-admissible. \\
a) If $k = 2$ then 
$$
\bar{d}\left(\{n \in \mb{N} : b_1(n+a_1) < b_2(n+a_2)\}\right) > 0,
$$
and in fact
$$
\limsup_{x \ra \infty} \frac{|\{n \leq x : n+a_1,n+a_2 \in \mc{N}, b_1(n+a_1) < b_2(n+a_2)\}|}{|\{n \leq x : n+a_1, n+a_2 \in \mc{N} \}|} \geq \frac{1}{2}.
$$
b) If $k = 3$ and at least two of $f_1,f_2,f_3$ are the same then
$$
\bar{d}\left(\{n \in \mb{N} : b_1(n+a_1) < b_2(n+a_2) < b_3(n+a_3)\} \right) > 0,
$$
and in fact
$$
\limsup_{x \ra \infty} \frac{|\{n \leq x : n+a_1,n+a_2,n+a_3 \in \mc{N}, b_1(n+a_1) < b_2(n+a_2) < b_3(n+a_3)\}|}{|\{n \leq x : n+a_1, n+a_2,n+a_3 \in \mc{N} \}|} \geq \frac{1}{6}.
$$
\end{thm}
An immediate corollary of this result is as follows.
\begin{cor} \label{cor:signedTau}
For any permutation $\sg\in S_3$ the sets 
$$
\{n \in \mb{N} : \tau(n + \sg(1)) < \tau(n+\sg(2)) < \tau(n+\sg(3))\}
$$
each have positive upper density. Assuming Lehmer's conjecture, this upper density is at least $1/6$.
\end{cor}
\begin{rem}
To the best of our knowledge this is the first non-trivial result towards Conjecture~\ref{conj:serre} for $k>2,$ where previously even the existence of infinitely many $n$ was not known. For the case $k=2,$ one could deduce that $b_1(n+a_1)<b_2(n+a_2)$ occurs with positive upper density simply by prescribing the signs $\text{sign}(b_1(n+a_1))=-1$ and $\text{sign}(b_2(n+a_2))=1$ using the work~\cite{TaTe} or by modifying somewhat an earlier approach from~\cite{MR}. This however would not allow us to get the right constant $1/2$ from the statement in a).
\end{rem}
The constraint $k \leq 3$ in Theorem \ref{thm:uncondSigns} results from our current limited understanding of correlations of bounded multiplicative functions. Assuming a conjecture of Elliott (generalizing a famous old conjecture of Chowla on correlations of the M\"{o}bius function), we may remove this constraint and also upgrade the conclusions from positive \emph{upper} density results to positive \emph{natural} density results. \\
Elliott's conjecture (as repaired in Section 1 of \cite{MRT}) can be phrased as follows. We denote below $\mb{U} := \{z \in \mb{C} : |z| \leq 1\}$.
\begin{conj}[\cite{MRT}, p. 5]\label{conj:EllIntro}
Let $g_1,\ldots,g_k:\mb{N} \ra \mb{U}$ be multiplicative functions and let $a_1,\ldots,a_k$ be distinct non-negative integers. Assume that there is an index $1 \leq j_0 \leq k$ such that for any fixed Dirichlet character $\chi$ we have
\begin{equation}\label{eq:nonPret}
\lim_{X \ra \infty} \inf_{|t| \leq X} \sum_{p \leq X} \frac{1-\text{Re}(g_{j_0}(p)\bar{\chi}(p)p^{-it})}{p} = \infty.
\end{equation}
Then as $X \ra \infty$,
$$
\sum_{n \leq X} g_1(n+a_1) \cdots g_k(n+a_k) = o(X).
$$
\end{conj}
We have the following result.
\begin{thm}\label{thm:condSigns}
Assume Conjecture \ref{conj:EllIntro}. Let $k\geq 2$, let $f_1,\ldots,f_k$ be arithmetically normalized primitive non-CM Hecke cusp forms and let $\mbf{a} = (a_1,\ldots,a_k)$ be a $\mbf{f}$-admissible tuple of distinct non-negative integers. Assume that the the hypotheses of Theorem \ref{thm:genCusp} are fulfilled. Then
$$
d\left(\{n \in \mb{N} : b_1(n+a_1) < \cdots < b_k(n+a_k)\}\right) > 0.
$$
In fact, as $x \ra \infty$ we have
$$
\lim_{x \ra \infty} \frac{|\{n \leq x : n+a_1,\ldots,n+a_k \in \mc{N}, b_1(n+a_1) < \cdots < b_k(n+a_k)\}|}{|\{n \leq x : n+a_1,\ldots,n+a_k \in \mc{N} \}|} = \frac{1}{k!}.
$$
\end{thm}

\subsection{Main technical propositions.}
\subsubsection{Unsigned Orderings} 
Theorem \ref{thm:genCusp} will follow from the following Erd\H{o}s-Kac type theorem. Below, we set $c := \frac{1}{2}(1+\pi^2/6)$. 
\begin{prop} \label{prop:EKCusp}
Let $k \geq 1$. Let $f_1,\ldots,f_k$ be primitive holomorphic non-CM cusp forms of weight $m$, squarefree level and trivial nebentypus, with respective sequences of Fourier coefficients $\{b_j(n)\}_n$, and assume that each of these satisfies assumption (A) from Theorem \ref{thm:genCusp}. Let $a_1,\ldots,a_k$ be distinct non-negative integers such that $\mbf{a}$ is $\mbf{f}$-admissible. Let $I_1,\ldots,I_k$ be closed non-empty intervals. Set $\mc{B}(\mbf{I}) := \prod_{1 \leq j \leq k} I_j$, and put $\lambda_j(n) := b_j(n)n^{-(m-1)/2}$ for all $n \in \mb{N}$ and $1 \leq j \leq k$. We write
$$
\mc{N}_j := \{n \in \mb{N} : \lambda_j(n) \neq 0\},
$$
and set $\mc{N} := \bigcap_{1 \leq j \leq k} \mc{N}_j$.
Finally, for each $1 \leq j \leq k$ and $n \in \mc{N}$ define
$$
X_j(n) := \frac{\log|\lambda_j(n)| + \frac{1}{2}\log\log n}{\sqrt{c\log\log n}}.
$$
Then, as $x \ra \infty$,
\begin{align*}
&\frac{1}{X}|\{n \leq x : n+a_j \in \mc{N}, X_j(n+a_j) \in I_j \text{ for all $1 \leq j \leq k$} \} | \\
&= (2\pi)^{-k/2}\int_{\mc{B}(\mbf{I})} e^{-\frac{1}{2}(x_1^2 + \cdots + x_k^2)} dx_1\cdots dx_k + o(1),
\end{align*}
where we have set
$$
X = X(x) := \sum_{n \leq x} \prod_{1 \leq j \leq k} 1_{\mc{N}}(n+a_j).
$$
\end{prop}
\begin{rem}
When $k = 1$ this proposition is due to Luca, Radziwi\l\l \ and Shparlinski \cite[Theorem 3]{LRS} (see also \cite{Ell} for a related result). More precisely, these authors proved this result assuming a certain quantitative form of the Sato-Tate conjecture that was subsequently derived by Thorner \cite{Tho}, based on the recent deep breakthrough of Newton-Thorne \cite{NeTh} on the functoriality of symmetric power $L$-functions for $\text{GL}(n)$ for all $n \geq 2$. Our result generalizes this to all fixed $k \geq 2$.
\end{rem}

\subsubsection{Signed Orderings}
Let $\mbf{\eps} \in \{-1,+1\}^k$ and $I_1,\ldots,I_k \subset \mb{R}$ be non-empty intervals. Keeping the notation from the previous subsection,
write $\sg_j(n) := \text{sign}(\lambda_j(n))$ whenever $\lambda_j(n) \neq 0$. Let us define the set
$$
\mc{S}_{\mbf{a}}(\mbf{\eps},\mbf{I}) := \{n \in \mb{N} : n+a_j \in \mc{N}, \sg_j(n+a_j) = \eps_j \text{ and } X_j(n+a_j) \in I_j \text{ for all } 1 \leq j \leq k\},
$$
and the counting function
$$
N_{\mbf{a}}(x; \mbf{\eps},\mbf{I}) := |\mc{S}_{\mbf{a}}(\mbf{\eps},\mbf{I}) \cap [1,x]|.
$$
We shall prove the following.
\begin{prop}\label{prop:EKwithSigns}
Let $\eps \in \{-1,+1\}^k$ and let $I_1,\ldots,I_k \subset \mb{R}$ be non-empty intervals. Write $\mc{B}(\mbf{I}) := \prod_{1 \leq j \leq k} I_j$ as above. \\
a) Assume Conjecture \ref{conj:EllIntro}. Then, provided $x$ is large enough in terms of $k,\mbf{a}$ and $\mbf{I}$,
$$
X^{-1}N_{\mbf{a}}(x;\mbf{\eps},\mbf{I}) = 2^{-k} \cdot (2\pi)^{-k/2} \int_{\mc{B}(\mbf{I})} e^{-\frac{1}{2}(x_1^2 + \cdots + x_k)^2} d\mbf{x} + o_{k,\mbf{I},\mbf{a}}(1).
$$ \\
b) Assume that $k \in \{2,3\}$ and that the hypotheses of Theorem \ref{thm:uncondSigns} hold (depending on $k$). Then 
$$
\sum_{n \leq x} \frac{1_{\mc{S}_{\mbf{a}}(\mbf{\eps},\mbf{I})}(n)}{n} = \left(2^{-k} \cdot (2\pi)^{-k/2}\int_{\mc{B}(\mbf{I})} e^{-\frac{1}{2}(x_1^2 + \cdots + x_k)^2} d\mbf{x} + o_{k,\mbf{I},\mbf{a}}(1) \right) \left(\sum_{n \leq x} \frac{1}{n}\prod_{1 \leq j \leq k} 1_{\mc{N}}(n+a_j)\right).
$$
\end{prop}
\begin{rem}
By the Sato-Tate theorem, the probability that $\sg_j(p) = +1$ or $-1$, for $p$ prime, is asymptotically $1/2$ in either case. Consequently, it is not hard to prove\footnote{This can be shown, e.g., when $\eps = +1$ by writing $1_{\sg_j(n) = \eps} = (1+\sg_j(n))/2$. The sum of $\sg_j(n)$ over $n \leq x$ is $o(x)$ by combining Sato-Tate with a result of Hall and Tenenbaum \cite{HaT}, for instance.} that the same is true of $\sg_j(n)$, for integers $n$ more generally. Note that, assuming that the events $\sg_j(n+a_j) = \e_j$ and $X_j(n+a_j) \in I_j$ are independent for each $j = 1,\ldots,k$, and moreover independent of the pairs $(\sg_{j'}(n+a_j),X_{j'}(n+a_{j'}))$ for all $j' \neq j$, then one would guess that the probability that $n \leq x$ satisfying $n+a_j \in \mc{N}_j$ for all $1 \leq j \leq k$ be counted by $N_{\mbf{a}}(x;\mbf{\eps},\mbf{I})$ is 
$$
2^{-k} (2\pi)^{-k/2} \int_{\mc{B}(\mbf{I})} e^{-\frac{1}{2}(x_1^2+\cdots + x_k^2)} d\mbf{x}.
$$
Thus, Proposition \ref{prop:EKwithSigns} a) (and b), in the case of logarithmic density) confirms this guess.
\end{rem}

\begin{rem}
The statistical behaviour of Fourier coefficients of Hecke cusp forms has been investigated extensively in various other aspects as well. For instance, in \cite{MMDuke} a variance estimate (in the spirit of the Hardy-Ramanujan theorem) was obtained for the number of distinct prime factors of non-zero prime and integer values of the Ramanujan tau function (the same authors extended these results further in \cite{MMEK}). Central limit theorems have also been deduced for partial sums of Fourier coefficients of cusp forms, see e.g., the work of Nagoshi \cite{Nag}, and the recent work of Murty and Prabhu \cite{MuPr}.
\end{rem}
\subsection{Organization of the paper and proof strategy} We shall now briefly describe the basic strategy underlying the proofs of our main results, focusing for simplicity on the case $b_f(n)=\tau(n)$. The starting point in our investigation is the simple observation that if $\tau(n) \neq 0$ then $\tau(n)=|\tau(n)|\text{sign}(\tau (n))$, where both $n \mapsto |\tau(n)|$ and $n\mapsto \text{sign}(\tau(n))$ are multiplicative functions. \\
In order to prove Proposition~\ref{prop:EKwithSigns}, after using standard Fourier analysis to approximate the indicator functions $1_{\tau(n+a_j)\in I_j},$ for $1\le j\le k$, roughly speaking we decouple the problem into understanding both the distribution of the vector 
$$
(\log|\tau(n+a_1)|,\log |\tau(n+a_2)|,\ldots,\log |\tau(n+a_k)|)
$$ 
as well as of the vector of signs $\text{sign}(\tau(n+a_i))=\epsilon_i$ where $\epsilon_i \in \{\pm 1\}$ for $i = 1,\ldots,k$.\\
To control the size of the components of the vector 
$$
(\log|\tau(n+a_1)|,\log |\tau(n+a_2)|,\dots,\log |\tau(n+a_k)|),
$$ 
we prove Proposition~\ref{prop:EKCusp}, showing that its distribution is jointly Gaussian (with trivial covariance matrix). In order to do so, we
first approximate the additive function $\log |\tau(n)|$ by the truncated version $\log |\tilde{\tau}_y(n)|$, which is supported on integers $n$ free of large prime factors $p>y$ and of those primes for which $\tau(p)$ is unusually small (see Lemma 4.1 for a more technical statement). This is accomplished by sieving out primes using the quantitative version of the Sato-Tate conjecture recently established by Thorner~\cite{Tho}, based on the breakthrough of Newton and Thorne~\cite{NeTh}. \\
In Section \ref{sec:aux}, we establish some consequences of the work~\cite{Tho} that will be used frequently in the remainder of the paper. In Section \ref{sec:sieve}, we use these results to establish sieve estimates to be used in the proof of Proposition \ref{prop:EKCusp}. In Section \ref{sec:mom}, we obtain an asymptotic formula for the moments of (an approximation of) the tuple $(X_1(n+a_1),\ldots,X_k(n+a_k))$, and deduce Proposition \ref{prop:EKCusp} using the method of moments in the manner of Granville and Soundararajan~\cite{GSEK}.\\
To study the frequency of the sign patterns 
$$
\text{sign}(\tau(n+a_1))=\epsilon_1,\text{sign}(\tau(n+a_2))=\epsilon_2,\dots,\text{sign}(\tau(n+a_k))=\epsilon_k,
$$
we use the recent, spectacular advances in the study of correlations of bounded multiplicative functions. It is here where, due to the limitations in our current understanding, we have to confine ourselves to sign patterns of length $k\le 3$ to get unconditional results. Specifically, to count the occurrences (with logarithmic weight) of the sign pattern $(\epsilon_1,\epsilon_2,\epsilon_3)$ we expand
\begin{align*}
\sum_{n\le x} \frac{\prod_{i\le 3}(1+\epsilon_i\text{sign}(\tau(n+a_i))}{n}&= (1+o(1))\log x+
\sum_{n\le x} \frac{\prod_{i\ne j\le 3}\epsilon_i\epsilon_j\text{sign}(\tau(n+a_i))\text{sign}(\tau(n+a_j))}{n}\\&+ \sum_{n\le x} \frac{\epsilon_1\epsilon_2\epsilon _3\text{sign}(\tau(n+a_1))\text{sign}(\tau(n+a_2))\text{sign}(\tau(n+a_3))}{n}
\end{align*}
and estimate\footnote{The logarithmic mean value terms are provably $o(\log x)$ by combining the Sato-Tate theorem with Hal\'{a}sz' theorem; see the proof of Proposition \ref{prop:EllApp} below.} the two- and three-point correlations of $\text{sign}(\tau(n))$ by using results of Tao~\cite{tao} and Tao-Ter\"av\"ainen~\cite{TaTe}, respectively. To this end we show that the function $\text{sign}(\tau(n))$ does not `` weakly pretend" to be $n^{it}\chi(n)$ for any Dirichlet character $\chi$, or more precisely
\[\min_{|t| \leq x} \sum_{p\le x}\frac{1-\Re{(p^{-it}\bar{\chi}(p)\text{sign}(\tau( p))})}{p}\gg \log \log x.\]
We do so by again appealing to the work of Newton and Thorne ~\cite{NeTh} and slightly modifying the results in~\cite{Tho} to handle the distribution of $\tau(p)$ in arithmetic progression $p\equiv a\pmod D$ for relatively small moduli $D$ (see Proposition~\ref{prop:ThoAPs}).\\
Our main applications, i.e., Theorem~\ref{thm:uncondSigns} and Theorem~\ref{thm:condSigns}, can straightforwardly be deduced from the above ideas, aside from obtaining the sharp density estimate $1/k!$ in each case. In Section 6, this is accomplished by a careful and somewhat involved argument of covering intervals with large Gaussian measure which might be of independent interest.
\begin{rem} It is possible to make further progress towards variants of Conjecture~\ref{conj:serre} for a more general class of multiplicative functions using methods developed here. More precisely, one can handle multiplicative functions $f:\mathbb{N}\to\mathbb{R},$ which can be represented as $f(n)=g(n)e^{h(n)}$ where $g$ is an oscillating multiplicative function that satisfies some form of Conjecture~\ref{conj:EllIntro} and $h$ is an additive function 
with Gaussian value distribution. In a similar fashion, one can prove related results restricting the shifts $n+a_1,\ldots, n+a_k$ to satisfy sieve theoretic conditions analogous to those considered in the work of Granville and Soundararajan \cite{GSEK}. In particular, as our proof shows, one can find a positive proportion of $n \leq x$, for which $P^+(n+a_j) \ll x^{\e}$ for all $1 \leq j \leq k$, that satisfy the conclusions of our main theorems. \\
We shall explore these lines of thought in a forthcoming work.
\end{rem}

\section*{Acknowledgments}
We would like to warmly thank Florian Luca and Maksym Radziwi\l\l \ for their helpful comments on an earlier draft of the present manuscript, and for their encouragement. The first author is grateful to the Max Planck Institute for Mathematics (Bonn) for providing excellent working conditions.

\section{Multiplicative functions supported on thin sets} \label{sec:aux}
Recall that in our setting, $f(z) := \sum_{n \geq 1} b_f(n)e(nz)$ is a primitive non-CM cusp form of weight $m$, square-free level $N$ and trivial nebentypus, and let $\lambda_f(n) := b_f(n)n^{-(m-1)/2}$. We define $\mc{N}_f := \{n \in \mb{N} : \lambda_f(n) \neq 0\}$, and set $\mc{N}_f(x) := \mc{N}_f \cap [1,x]$ for $x \geq 2$. \\
Since $\lambda_f(p) \in [-2,2]$ by deep work of Deligne \cite{Del}, we may define $\theta_f(p) := \cos^{-1}(\lambda_f(p)/2) \in [0,\pi]$. Given a closed interval $I \subseteq [0,\pi]$, we define
$$
\pi_{f}(x;I) := |\{p \leq x : \theta_f(p) \in I\}|,
$$
and the Sato-Tate measure
$$
\mu_{ST}(I) := \frac{2}{\pi} \int_I \sin^2 \theta d\theta.
$$
In the sequel, we will frequently have to deal with multiplicative functions supported on primes $p$ for which $|\lambda_f(p)| < (\log\log x)^{-1}$. The goal of this section is to develop estimates for partial sums of such functions.\\
The following result of Thorner, proved using the recent, deep breakthrough of Newton and Thorne \cite{NeTh} on functoriality for the symmetric power $L$-functions for $GL(n)$ and all $n \geq 2$, will be used crucially.
\begin{thm}[\cite{Tho}, Theorem 1.1] \label{thm:effST}
Let $I \subseteq [0,\pi]$ be a closed interval. Then
$$
\pi_{f}(x;I) = \left(\mu_{ST}(I)  + O\left(\frac{\log(kN\log x)}{(\log x)^{1/2}}\right)\right)\pi(x).
$$
\end{thm}
We next present a sharpening of a result from \cite{LRS}, which shows that the set of primes with $\lambda_f(p) = 0$ is quite thin.
\begin{lem} \label{lem:LRSSharp}
Let $y \geq 2$. 
Then
$$
\sum_{p > y \atop \lambda_f(p) = 0} \frac{1}{p} \ll_f \frac{(\log\log y)^2}{\sqrt{\log y}}.
$$
Also, for any $A \geq 1$ and $x \geq y$ large enough we have 
$$
\sum_{y < p \leq x \atop |\lambda_f(p)| < (\log\log x)^{-A}} \frac{1}{p} \ll_f (\log\log x)^{1-A}.
$$
\end{lem}
\begin{proof}
We begin with the first statement.
We may assume that $y$ is larger than any fixed constant depending at most on $f$ (by selecting a suitably large implicit constant). Let $\psi: \mb{R} \ra [0,\infty)$ be a non-increasing function, satisfying $\psi(y) \ra 0$ as $y \ra \infty$, to be chosen later. 
We clearly have
$$
\sum_{p > y \atop \lambda_f(p) = 0} \frac{1}{p} \leq \sum_{p > y \atop |\lambda_f(p)| \leq \psi(p)} \frac{1}{p}.
$$
We apply Theorem \ref{thm:effST} and partial summation. Write 
$$
E_{f}(X;I) := \pi_{f}(X;I) - \mu_{ST}(I)\pi(X),
$$
for any closed interval $I \subseteq [0,\pi]$. Set $k_0 := \llf \log \log y\rrf$, and $x_k := e^{e^k}$ for all $k \geq k_0$. For each $k \geq k_0$ put $J_k := \cos^{-1}([-\psi(x_k),\psi(x_k)]) \subseteq [0,\pi]$, noting that $|J_k| \asymp \psi(x_k)$ and $[-\psi(p),\psi(p)] \subseteq J_k$ whenever $x_k < p \leq x_{k+1}$. Then by partial summation,
\begin{align*}
\sum_{p > y \atop |\lambda(p)| \leq \psi(p)} \frac{1}{p} &\leq \sum_{k\geq k_0} \left(\mu_{ST}(J_k)\int_{x_k}^{x_{k+1}} \frac{du}{u\log u} + \int_{x_k}^{x_{k+1}} u^{-2}E_f(u;J_k) du + O\left(e^{-k}\right)\right) \\
&= \sum_{k \geq k_0} \mu_{ST}(J_k) + O\left(\frac{\log\log y}{\sqrt{\log y}}\right).
\end{align*}
For fixed $\e > 0$ we select $\psi(y) := (\log \log y)^2(\log y)^{-1/2}$, so that $\mu_{ST}(J_k) \asymp |J_k| \asymp k^2e^{-k/2}$. We thus find that
$$
\sum_{p > y \atop |\lambda_f(p)| \leq \psi(p)} \frac{1}{p} \ll k_0^2e^{-k_0/2} \ll \frac{(\log\log y)^2}{\sqrt{\log y}},
$$
as claimed. \\
The second claim follows similarly from Theorem \ref{thm:effST} by partial summation, given that
$$
\sum_{p \leq x \atop |\lambda_f(p)| < (\log\log x)^{-A}} \frac{1}{p} \sim \mu_{ST}([-(\log\log x)^{-A},(\log\log x)^{-A}]) \log\log x \ll (\log\log x)^{1-A},
$$
whenever $A \geq 1$.
\end{proof}

%

The next two lemmas are the main technical results of this section.
\begin{lem} \label{lem:thinSetProp}
Let $S$ be a set of primes such that $\sum_{p \in S} 1/p < \infty$. Let $g: \mb{N} \ra \mb{R}$ be a multiplicative function such that $g(p) = 0$ unless $p \in S$, and there is $k \geq 1$ for which $|g(n)| \leq d_k(n)$ for all $n$. The following holds:
\begin{itemize}
\item we have $\sum_{n \geq 1} \frac{|g(n)|}{n} < \infty;$ \\
\item for any $10 \leq y \leq z \leq x$,  put $u := \frac{\log z}{\log y}$. Then
$$\sum_{z < n \leq x} \frac{|g(n)|}{n} \ll_k \exp\left(-\frac{\log z}{\log y}\right) + \sum_{y < p \leq x \atop p \in S} \frac{1}{p} + y^{-1}.
$$
\end{itemize}
\end{lem}
\begin{proof}
The first statement easily follows from the uniform bound 
$$
\sum_{n \leq X} \frac{|g(n)|}{n} \leq \sum_{P^+(n) \leq X} \frac{|g(n)|}{n} \ll_k \exp\left(k\sum_{p \leq X \atop p \in \mc{S}} \frac{1}{p}\right)\ll_k 1.
$$
 To prove the second claim, we split the sum as
$$
\sum_{z < n \leq x} \frac{|g(n)|}{n} \leq \sum_{n > z \atop P^+(n) \leq y} \frac{|g(n)|}{n} + \sum_{z < n \leq x \atop P^+(n) > y} \frac{|g(n)|}{n}.
$$
Applying Rankin's trick with $\sg := 1/\log y < 1/2$ yields
$$
\sum_{n > z \atop P^+(n) \leq y} \frac{|g(n)|}{n} \leq z^{-\sg} \sum_{P^+(n) \leq y} \frac{|g(n)|}{n^{1-\sg}} \ll_k z^{-\sg} \exp\left(k\sum_{p \leq y \atop p \in S} \frac{p^{\sg}}{p}\right) \ll_k e^{-u}\exp\left(ky^{\sg}\right) \ll_k e^{-u}.
$$
For the second term, we note that
$$
\sum_{z < n \leq x \atop P^+(n) > y} \frac{|g(n)|}{n} \leq \sum_{y < p \leq x \atop \nu \geq 1} \frac{|g(p^{\nu})|}{p^{\nu}} \sum_{d\geq 1} \frac{|g(d)|}{d}.
$$
The inner sum is bounded by i), while the outer sum is, by partial summation and the bound $d_k(p^{\nu}) \ll_k \nu^k$,
\begin{align*}
&\ll_k \sum_{y < p \leq x \atop p \in S} \frac{1}{p} + \sum_{p > y \atop \nu \geq 2} \frac{d_k(p^{\nu})}{p^{\nu}} \\
&\ll_k \sum_{y < p \leq x \atop p \in S} \frac{1}{p} + \frac{1}{(\log y)^k}\sum_{p > y \atop \nu \geq 2} \frac{(\log(p^{\nu}))^k}{p^{\nu}} \ll \sum_{y < p \leq x \atop p \in S} \frac{1}{p} + y^{-1},
\end{align*}
and the claim follows.
\end{proof}
Let $1 \leq j \leq k$. Upper bounds for the partial sums $\sum_{n \leq x} |g_j(n)|$, where $g_j$ is supported on primes $p$ for which $|\lambda_j(p)| < (\log\log x)^{-A}$ as $x \ra \infty$, will play an important role in our arguments later on. We note that by Shiu's theorem and Lemma \ref{lem:LRSSharp}, we have the cheap bound
\begin{equation}\label{eq:cheap}
\sum_{n \leq x} |g_j(n)| \ll \frac{x}{\log x} \exp\left(\sum_{p \leq x \atop |\lambda_j(p)| < (\log\log x)^{-A}} \frac{1}{p}\right) \ll \frac{x}{\log x},
\end{equation}
for each $1 \leq j \leq k$. The next lemma provides an improvement to this bound.

\begin{lem}\label{lem:STThin}
Let $x$ be large and let $A \geq 1$. Let $\{\lambda_f(n)\}_n$ be the set of normalized Fourier coefficients of a primitive non-CM holomorphic cusp form. Let $h: \mb{N} \ra [0,1]$ be a multiplicative function such that 
$$
h(n) = 0 \text{ unless } |\lambda_f(p)| < (\log\log x)^{-A} \text{ whenever } p|n.
$$
Then
$$
\sum_{n \leq x} h(n) \ll_{f,A} \frac{x}{(\log x) (\log \log x)^A}.
$$
\end{lem}

\begin{proof}
We write $\Lambda_h$ to denote the coefficients of the Dirichlet series $-L'(s,h)/L(s,h)$, where 
$$
L(s,h) := \sum_{n \geq 1} h(n)n^{-s} \text{ for $\text{Re}(s) > 1$.}
$$
We note that $\Lambda_h(p) = h(p)\log p$ for all primes $p$, and that $|\Lambda_h(n)| \leq \Lambda(n)$ for all $n.$
Applying the hyperbola method, we obtain
\begin{align} \label{eq:logSmooth}
\sum_{n \leq x} h(n)\log n = \sum_{a \leq \sqrt{x}} h(a) \sum_{b \leq x/a} \Lambda_h(b) + \sum_{b \leq \sqrt{x}} \Lambda_h(b) \sum_{\sqrt{x} < a \leq x/b} h(b) =: T_1 + T_2.
\end{align}
To proceed, we note that by partial summation and Theorem \ref{thm:effST}, we have, for $\sqrt{x} \leq y \leq x$,
\begin{align*}
\sum_{m \leq y} \frac{\Lambda_h(m)}{m} = \sum_{p \leq y} \frac{h(p)\log p}{p} + O_f(1) \leq \sum_{p \leq y \atop |\lambda_f(p)| < (\log\log x)^{-A}} \frac{\log p}{p} + O_f(1) \ll_{A,f} \frac{\log y}{(\log\log x)^A},
\end{align*}
as well as 
\begin{align*}
\sum_{m \leq y} \Lambda_h(m) = \sum_{p \leq y \atop |\lambda_f(p)| < (\log \log x)^{-A}} h(p)\log p +O_f(\sqrt{y} \log y) 
\ll_{A,f} \frac{y}{(\log\log x)^A}.
\end{align*}
Applying these estimates, \eqref{eq:cheap} and Lemma \ref{lem:LRSSharp} to \eqref{eq:logSmooth}, we obtain
\begin{align*}
\sum_{n \leq x} h(n)\log n &\ll x\sum_{a \leq \sqrt{x}} \frac{h(a)}{a(\log\log x)^A} + x\sum_{b \leq \sqrt{x}} \frac{\Lambda_h(b)}{b\log x} \\
&\ll \frac{x}{(\log\log x)^A} \exp\left(\sum_{p \leq x} \frac{h(p)}{p}\right) + \frac{x}{\log x} \cdot \frac{\log x}{(\log\log x)^A} \\
&\ll_A \frac{x}{(\log\log x)^A}.
\end{align*}
It follows that
$$
\sum_{n \leq x} h(n) \leq \sum_{n \leq \sqrt{x}} h(n)\log n + \frac{2}{\log x} \sum_{\sqrt{x} < n \leq x} h(n)\log n \ll_{A,f} \frac{x}{(\log x)(\log\log x)^A},
$$
as claimed.
\end{proof}

%
\section{Sieve Estimates for Support Sets of Coefficients of Cusp Forms} \label{sec:sieve}
We recall our assumption (A) from the statement of Theorem \ref{thm:genCusp}, from which we know that the following auxiliary assumption also holds: 
$$
\text{(A'): } p \in \mc{N}_j \Rightarrow p^k \in \mc{N}_j \text{ for all $1 \leq j \leq k$}.
$$ 
This condition also clearly holds with $\mc{N} = \bigcap_{1 \leq j \leq k} \mc{N}_j$ in place of the $\mc{N}_j$.
In addition, assumption (A) requires that $|\lambda_j(p^k)| \geq p^{-kC}$, for some $C = C(f_j) > 0$, a condition that will be used in the next section. \\
Given a non-increasing function $\xi: \mb{N} \ra \mb{R}$ with $\xi > 0$ and $x$ large we also define 
$$
\mc{N}^{(\xi)}(x) := \bigcap_{1 \leq j \leq k} \{n \leq x : n\in \mc{N}_j,  p|n \Rightarrow |\lambda_j(p)| \geq \xi(x)\}.
$$
Note that $\mc{N}^{(\xi)}(x) \subseteq \mc{N} \cap [1,x]$, and that $\mc{N}^{(\xi)}(x)$ is a multiplicative set, i.e., if $m,n \in \mb{N}$ are coprime then $mn \in \mc{N}^{(\xi)}(x)$ if and only if $m,n \in \mc{N}^{(\xi)}(x)$.  In typical applications below we will choose $\xi(x) = (\log\log x)^{-A}$, for some $A > 1$. \\
Lastly, in what follows we set
$$
\Delta_{\mbf{a}} := \prod_{1 \leq i < j \leq k} (a_j-a_i),
$$
which is the discriminant of the system of linear forms $\{n \mapsto n+a_j\}_{1 \leq j \leq k}$. \\
The main goal of this section is to prove the following estimate, which will play a crucial role in our moment computations in the next section.
\begin{prop} \label{prop:Sieve}
Let $x$ be large. Assume that $\mbf{a} = (a_1,\ldots,a_k)$ is $\mbf{f}$-admissible. Then there are multiplicative functions $F_j :\mb{N} \ra [0,\infty)$ such that the following holds: given $d_1,\ldots, d_k \in \mb{N}$ coprime to $\Delta_{\mbf{a}}$ with $d_j \in \mc{N}^{(\xi)}(2x)$ and $D := d_1\cdots d_k \leq x^{1/4}$, then
$$
\sum_{\substack{n \leq x \\ d_j|(n+a_j) \\ \forall \ 1 \leq j \leq k}} \prod_{1 \leq j \leq k} 1_{\mc{N}^{(\xi)}(2x)}(n+a_j) = X\frac{F_1(d_1)\cdots F_k(d_k)}{d_1\cdots d_k} + O(r_{D}(x)),
$$
where $r_m(x) \ll \frac{x}{m(\log \log x)^{A-1}}$ for $m \geq 1$, and 
$$
X = X(x) := \sum_{n \leq x} \prod_{1 \leq j \leq k} 1_{\mc{N}}(n+a_j) \gg_{\mbf{a},k} x.
$$
\end{prop}

We need the following generalization of Lemma 8 of \cite{LRS}.
\begin{lem}\label{lem:approx}
Let $I_1,\ldots,I_k$ be a sequence of intervals. Assume that $\xi(x) \geq (\log x)^{-1/2+\e}$, for any $\e > 0$. Then
\begin{align*}
&|\{n \leq x : n+a_j \in \mc{N}, \log|\lambda_{j}(n+a_j)| \in I_j \ \forall \ 1 \leq j \leq k\}| \\
&= |\{n \leq x : n+a_j \in \mc{N}^{(\xi)}(2x), \log|\lambda_{j}(n+a_j)| \in I_j \ \forall \ 1 \leq j \leq k\}| + O_{k,\mbf{f}}\left(x(\log\log x) \xi(x)\right).
\end{align*}
\end{lem}
\begin{proof}
We denote by $N(x)$ and $N_{\xi}(x)$ the cardinalities on the LHS and RHS, respectively. 
As $\mc{N}^{(\xi)}(2x) \subseteq \mc{N}$, we have $N_{\xi}(x) \leq N(x)$. By the union bound, 
$$
0 \leq N(x) - N_{\xi}(x) \leq \sum_{1 \leq j \leq k} |\{n \leq x : n+a_j \in \mc{N} \bk \mc{N}^{(\xi)}(2x)\}|.
$$
The latter cardinality counts integers $n$ having at least one prime factor $p$ such that $0 < |\lambda_j(p)| \leq \xi(x)$, for some $1 \leq j \leq k$. Partial summation as in Lemma \ref{lem:LRSSharp} yields
\begin{align*}
|\{n \leq x : n+a_j \in \mc{N} \bk \mc{N}^{(\xi)}(2x)\}| &\leq 2x\sum_{p \leq 2x \atop \lambda_j(p) \in [-\xi(x),\xi(x)]} \frac{1}{p} \\
&\ll_{\mbf{f}} x(\log\log x)\left(\xi(x) + \frac{\log\log x}{\sqrt{\log x}}\right) \\
&\ll x(\log\log x) \xi(x).
\end{align*}
This implies the claim.
\end{proof}

We fix $A > 1$ (taken to be larger later), and set $\xi(y) := (\log\log y)^{-A}$.
\begin{lem} \label{lem:badDiv}
Let $x$ be large, and let $10 \leq z \leq x$. Define $\tilde{g} := \mu \ast 1_{\mc{N}^{(\xi)}(2x)}$. Let $d_1,\ldots,d_k$ be integers coprime to $\Delta_{\mbf{a}}$, $d_j \in \mc{N}^{(\xi)}(2x)$ for all $j$, and such that $D:= d_1\cdots d_k \leq x^{1/2}$. Then the number of $n \leq x$ such that:
\begin{itemize}
\item $d_j|n+a_j$ for all $1 \leq j \leq k$, and
\item there is $j_0 \in \{1, \ldots, k\}$ and $e > z$ such that $e|(n+a_{j_0})$ and $\tilde{g}_{j_0}(e) \neq 0$
\end{itemize} 
satisfies the upper bound
$$
\ll \frac{x}{D}\left( (\log z)^{-1/2+o(1)} + (\log\log x)^{1-A} \right).
$$
\end{lem}
\begin{proof}
We first observe that
$$
(d_i,d_j) = (e_i,d_j) = (e_i,d_i) = 1 \text{ for all $i \neq j$,} 
$$ 
whenever $e_i|(n+a_i)$ and $\tilde{g}(e_i) \neq 0$. Indeed, the first equality arises from the fact that $(d_i,d_j)$ divides $(n+a_i)-(n+a_j) = a_i-a_j$, which is a divisor of $\Delta_{\mbf{a}}$. The second equality arises for the same reason. To verify the third equality note that if $p|(e_i,d_i)$ then $p \in \mc{N}_i^{(\xi)}(2x)$, yet $\tilde{g}(e_i) = \prod_{p'||e_i}(1-1_{\mc{N}^{(\xi)}(2x)}(p'))$ vanishes. The conclusion from this is that whenever there is $1 \leq i \leq k$ and a divisor $e_i|n+a_i$ with $\tilde{g}(e_i) \neq 0$ then $(e_i,[d_1,\ldots,d_k]) = (e_i,d_1\cdots d_k) = 1$. \\
By the Chinese remainder theorem there is a residue class $b \pmod{D}$ such that the cardinality we seek to estimate equals
\begin{align*}
|\{n \leq x : n \equiv b \pmod{D} \text{ and } \exists \ j \in \{1,\ldots,k\}, &e_j > z, (e_j,D) = 1 \text{ such that } \\
&\tilde{g}(e_j) \neq 0 \text{ and } e_j|(n+a_j) \}|,
\end{align*}
and by the union bound, this is 
$$
\leq \sum_{1 \leq j \leq k} |\{n \leq x : n \equiv b \pmod{D} \text{ and } \exists \ e > z, (e,D) = 1 \text{ such that } \tilde{g}(e) \neq 0 \text{ and } e|(n+a_j)\}|.
$$
Fix $1 \leq j \leq k$. Note that if $\tilde{g}(e) \neq 0$ then $|\tilde{g}(e)| = 1$. The $j$th term on the RHS in this last inequality is, again by the union bound,
\begin{align*}
&\leq \sum_{z < e \leq x/D \atop (e,D) = 1} |\tilde{g}(e)| \sum_{\substack{n \leq x \\ e|(n+a_j) \\ D|(n-b)}} 1 + \sum_{x/D < e \leq 2x \atop (e,D) = 1} |\tilde{g}(e)| \sum_{\substack{n \leq x \\ e|(n+a_j) \\ D|(n-b)}} 1 =: T_{1,j} + T_{2,j}.
\end{align*}
We first estimate the sum over $T_{1,j}$. Using the Chinese remainder theorem, when $e \leq x/D$ is coprime to $D$ the inner sum is $\ll x/(eD)$. Now set $y := \exp((\log z)/(\log\log z))$. Since $|\tilde{g}(p)|$ is supported on $p \notin \mc{N}^{(\xi)}(2x)$, we may apply Lemma \ref{lem:thinSetProp}ii) and Lemma \ref{lem:LRSSharp} to get
\begin{align*}
\sum_{1 \leq j \leq k} T_{1,j} &\ll_k \frac{x}{D} \sum_{z < e \leq x} \frac{|\tilde{g}(e)|}{e} \ll_k \frac{x}{D}\left(\exp\left(-\frac{\log z}{\log y}\right) + y^{-1} + \sum_{1 \leq j \leq k} \left(\sum_{p > y \atop p \notin \mc{N}_j} \frac{1}{p} + \sum_{y < p \leq x \atop |\lambda_j(p)| <\xi(2x)} \frac{1}{p}\right) \right) \\
&\ll \frac{x}{D}\left( (\log y)^{-1/2+o(1)} + (\log\log x)^{1-A}\right) \\
&\ll \frac{x}{D}\left( (\log z)^{-1/2+o(1)} + (\log\log x)^{1-A}\right),
\end{align*}
which is sufficient. Next, consider the terms $T_{2,j}$, which we can bound as
\begin{align*}
T_{2,j} \leq \sum_{x/D < e \leq 2x \atop (e,D) = 1} |\tilde{g}(e)| \sum_{m \leq x/D \atop e|mD+b+a_j} 1 = \sum_{x/D < e \leq 2x \atop (e,D) = 1} \frac{|\tilde{g}(e)|}{e} \sum_{0 \leq r \leq e-1} e\left(\frac{r(b+a_j)}{e}\right) \sum_{m \leq x/D} e\left(\frac{rmD}{e}\right),
\end{align*}
for each $1 \leq j \leq k$. Separating the term $r = 0$ from the rest and then using a geometric series bound, we obtain
\begin{align*}
\sum_{1 \leq j \leq k} T_{2,j} &\ll_k \frac{x}{D}\sum_{x/D < e \leq 2x} \frac{|\tilde{g}(e)|}{e} + \sum_{x/D < e \leq 2x \atop (e,D) = 1} \frac{|\tilde{g}(e)|}{e} \sum_{1 \leq r \leq e-1} \min\{x/D,\|rD/e\|^{-1}\} \\
&\ll \frac{x}{D}(\log\log x)^{1-A} + \frac{\log x}{D}\sum_{x/D < e \leq 2x} |\tilde{g}(e)|,
\end{align*}
where we bounded the first term in the last line as we did for $T_{1,j}$ (but with $x/D \geq x^{1/2}$ in place of $z$). To treat the second term we apply Lemma \ref{lem:STThin}, since $|\tilde{g}|$ is supported only on primes $p$ with $|\lambda_j(p)| < (\log\log x)^{-A}$ for some $1 \leq j \leq k$. This gives
$$
\sum_{1 \leq j \leq k} T_{2,j} \ll_k \frac{x}{D} (\log\log x)^{1-A}.
$$
The claim follows upon combining the bounds for $T_{1,j}$ and $T_{2,j}$.
\end{proof}
%
We are now ready to complete the proof of the main proposition of this section.

\begin{proof}[Proof of Proposition~\ref{prop:Sieve}]
Let $1 \leq y \leq x^{1/(4k)}$. Lemma \ref{lem:badDiv} implies that
$$
\mc{S}_{\mbf{d}}(x) := \sum_{\substack{n \leq x \\ d_j |(n+a_j) \forall j \\ 1\leq j \leq k}} \prod_{1 \leq j \leq k} 1_{\mc{N}^{(\xi)}(2x)}(n+a_j) = \asum_{\substack{n \leq x \\ d_j |(n+a_j) \\ 1\leq j \leq k}} \prod_{1 \leq j \leq k} 1_{\mc{N}^{(\xi)}(2x)}(n+a_j) + O_k\left(\frac{x}{D(\log \log x)^{A-1}}\right),
%
$$
where the asterisk in the sum means that $n$ has the additional constraint that if, for any $1\leq j \leq k$, $e | (n+a_j)$ and $\tilde{g}(e) \neq 0$ then $e \leq y$. 
Using the definition of $\tilde{g}$ we may write
$$
\mc{S}_{\mbf{d}}(x) = \sum_{\substack{e_1,\ldots,e_k \leq y \\ (e_id_i,e_jd_j) = 1 \\ \forall i \neq j}} \tilde{g}(e_1)\cdots \tilde{g}(e_k) \asum_{\substack{n \leq x \\ e_jd_j | n+a_j \\ 1\leq j \leq k}} 1 + O_k\left(\frac{x}{D(\log \log x)^{A-1}}\right),
$$
where, as in the previous lemma, we note that $(e_i,d_j) = 1$ for all $1 \leq i,j \leq k$, and moreover $(e_i,e_j) = 1$ for all $i \neq j$, for in all other cases the summands vanish. According to assumption (A), $\tilde{g}$ is supported on squarefree integers, and $|\tilde{g}(m)| \in \{0,1\}$ for all $m \in \mb{N}$. 
By the Chinese remainder theorem, for each tuple $(e_1,\ldots,e_k)$ there is a residue $c \pmod{d_1e_1\cdots d_ke_k}$ such that $n \equiv c \pmod{d_1e_1 \cdots d_ke_k}$, noting that $e_1d_1 \cdots e_kd_k \leq x^{1/2}$ by choice of $y$. We conclude, by a second application of Lemma \ref{lem:badDiv} (with $e_jd_j$ in place of $d_j$ for all $j$) that
\begin{align*}
\mc{S}_{\mbf{d}}(x) &= \sum_{\substack{e_1,\ldots,e_k \leq y \\ (e_id_i,e_jd_j) = 1 \\ \forall i \neq j}} \tilde{g}(e_1)\cdots \tilde{g}(e_k) \frac{x}{e_1\cdots e_kD} \left(1+O_k\left(\frac{1}{(\log\log x)^{A-1}}\right)\right) + O_k\left(\frac{x}{D(\log \log x)^{A-1}}\right)\\
&= \frac{x}{D}\sum_{\substack{e_1,\ldots,e_k \leq y\\ (e_id_i,e_jd_j) = 1 \forall i \neq j}} \frac{\tilde{g}(e_1)\cdots \tilde{g}(e_k)}{e_1\cdots e_k} + O_k\left(\frac{x}{D(\log\log x)^{A-1}}\left(1+ \prod_{1 \leq j \leq k} \sum_{e_j \leq y} \frac{|\tilde{g}(e_j)|}{e_j}\right)\right).
\end{align*}
We may bound the error term easily, via Lemma \ref{lem:LRSSharp}, giving
$$
\ll \frac{x}{D} (\log\log x)^{1-A}\exp\left(\sum_{1 \leq j \leq k} \sum_{p \leq y} \frac{|\tilde{g}(p)|}{p}\right) \ll \frac{x}{D}(\log\log x)^{1-A}.
$$
%
We can express the main term as
$$
\frac{x}{D} \left(\sum_{\substack{e_1,\ldots,e_k \leq y \\ (e_i,e_j) = 1 \forall i \neq j \\ (e_i,D) = 1 \forall i}} \frac{\tilde{g}(e_1)\cdots \tilde{g}(e_k)}{e_1\cdots e_k}\right) = \frac{x}{D}\left(\sum_{\substack{e_1,\ldots,e_k \leq y \\ (e_i,e_j) = 1 \forall i \neq j}} \frac{\tilde{g}(e_1)\cdots \tilde{g}(e_k)}{e_1\cdots e_k}\right).
$$
the condition $(e_i,D) = 1$ being redundant in view of $D \in \mc{N} \cap [1,x]$ and $\tilde{g}$ vanishing on that set. Thus,
%
%
$$
S_{\mbf{d}}(x) = \frac{x}{d_1\cdots d_k} \left(\sum_{e_1,\ldots,e_k \leq y \atop (e_i,e_j) = 1 \forall i \neq j} \frac{\tilde{g}(e_1)\cdots \tilde{g}(e_k)}{e_1\cdots e_k} + O_k\left((\log \log x)^{1-A}\right)\right).
$$
Applying the same argument when $d_j = 1$ for all $1 \leq j \leq k$, followed by Lemma \ref{lem:approx} (with $I_j = \mb{R}$ for all $j$), we find
\begin{align*}
x\sum_{\substack{e_1,\ldots,e_k \leq y \\ (e_i,e_j) =1 \\ \forall i \neq j}} \frac{\tilde{g}(e_1)\cdots \tilde{g}(e_k)}{e_1\cdots e_k} &=  S_{(1,\ldots,1)}(x) + O_k(x(\log\log x)^{1-A}) \\
&= \sum_{n \leq x} \prod_{1 \leq j \leq k}1_{\mc{N}^{(\xi)}(2x)}(n+a_j) + O_k(x(\log\log x)^{1-A}) \\
&= X+ O_k(x(\log\log x)^{1-A}), 
\end{align*}
and the claimed estimate follows with $F_j(d) := 1_{\mc{N}^{(\xi)}(2x)}(d)$ for each $1 \leq j \leq k$ (that $F_j(d_j) = 1$ for all $1 \leq j \leq k$ is trivially satisfied given the hypotheses). \\
To verify that $X \gg_{\mbf{a},k} x$ we argue differently. By Lemma \ref{lem:LRSSharp}, we find that
$$
\sum_{\log x < p \leq x \atop p \notin \mc{N}} \frac{1-1_{\mc{N}^{(\xi)}(2x)}(p)}{p} \leq \sum_{ 1 \leq j \leq k} \sum_{ \log x <p \leq x \atop p \notin \mc{N}} \frac{1_{|\lambda_j(p)| < (\log\log x)^{-A}}}{p} \ll \frac{(\log\log \log x)^2}{(\log\log x)^{1/2}},
$$
for each $1 \leq j \leq k$. By Theorem 1.3 of \cite{klu} (or, more precisely, the more general Theorem 2.1 of \cite{KM}), we have
$$
X = \sum_{n \leq x} \prod_{1 \leq j \leq k} 1_{\mc{N}}(n+a_j) = x\left(\prod_{p \leq x} M_p(\mbf{a}) + O_k\left(\frac{(\log \log \log x)^2}{(\log\log x)^{1/2}} \right)\right),
$$
where, writing $\nu_p(n) := \max\{ \nu \geq 0 : p^{\nu}|n\}$ for each prime $p$ and $n \in \mb{N}$, we have set
$$
M_p(\mbf{a}) := \lim_{z \ra \infty} z^{-1} \sum_{n \leq z} \prod_{1 \leq j \leq k} 1_{\mc{N}}(p^{\nu_p(n+a_j)}).
$$
It thus suffices to show that
$$
\prod_{p \leq x} M_p(\mbf{a}) \gg_{\mbf{a},k} 1.
$$
For $p \nmid \Delta_{\mbf{a}}$ with $p \in \mc{N}$ we easily find
$$
M_p(\mbf{a}) = 1-\frac{k}{p} + \sum_{1 \leq j \leq k} \lim_{z \ra \infty} \sum_{\nu \geq 1} z^{-1}\sum_{n \leq z \atop p^{\nu}||n+a_j} 1 = 1-\frac{k}{p} + k \left(1-\frac{1}{p}\right) \sum_{j \geq 1} \frac{1}{p^j} = 1,
$$
so that as $1_{\mc{N}_j}(p) > 0$ for all $j$ we must have 
\begin{equation}\label{eq:goodProd}
\prod_{\substack{p \leq x \\ p \in \mc{N}}} M_p(\mbf{a}) \gg_{\mbf{a}} \prod_{\substack{p \leq x \\ p \in \mc{N}\\ p \nmid \Delta_{\mbf{a}}}} M_p(\mbf{a}) \gg 1.
\end{equation}
Furthermore, since $\mbf{a}$ is $\mbf{f}$-admissible, when $p \notin \mc{N}$ we have
$$
M_p(\mbf{a}) = \lim_{z \ra \infty} z^{-1}\sum_{\substack{n \leq z \\ p \nmid \prod_j(n+a_j)}} 1 \gg \frac{1}{p}.
$$
Therefore, as
$$
\sum_{p \leq x \atop p \notin \mc{N}} \frac{1}{p} 
\leq \sum_{1 \leq j \leq k} \sum_{p \leq x \atop p \notin \mc{N}_j} \frac{1}{p} \ll_{\mbf{a},k} 1
$$ 
by Lemma \ref{lem:LRSSharp}, we deduce that
\begin{align*}
\prod_{p \leq x \atop p \notin \mc{N}} M_p(\mbf{a}) \gg \exp\left(-\sum_{1 \leq j \leq k} \sum_{p \leq x \atop p \notin \mc{N}_j} \frac{1}{p}\right) \gg_k 1.
\end{align*}
Together with \eqref{eq:goodProd}, this implies the bound $X \gg_{\mbf{a},k} x$, as required.
\end{proof}
\section{Moment Computation and Proof of Proposition \ref{prop:EKCusp}} \label{sec:mom}
In light of assumption (A), for $y \geq 2$ we define for each $n \in \mc{N}$
$$
\tilde{\lambda}_{j,y}(n) := \prod_{\substack{p|n \\ p \leq y \\ p \nmid \Delta_{\mbf{a}}}} \lambda_j(p).
$$ 
Throughout this section, let $A > 1$ be large and let $\xi(x) := (\log\log x)^{-A}$. We also assume henceforth that $\mbf{a}$ is $\mbf{f}$-admissible. We write $\Delta$ in place of $\Delta_{\mbf{a}}$ for convenience. The next lemma shows that $\tilde{\lambda}_{j,y}(n)$ well-approximates $\lambda_j(n)$ for most  choices $n\le x.$
\begin{lem} \label{lem:assump}
Let $2 \leq y \leq x$ with $y \geq x^{1/\log\log x}$. For all but $o_{\mbf{a},\mbf{f},A}(x)$ integers $n \in \mc{N}^{(\xi)}(x)$ we have 
$$
\log|\lambda_j(n)| = \log|\tilde{\lambda}_{j,y}(n)| + O( (\log\log\log x)^3) \text{ for all $1 \leq j \leq k$.}
$$
\end{lem}
\begin{proof}
We use an idea from \cite[Lemma 9]{LRS}. If $n \in \mc{N}$ and $p^{\nu}|| n$ then by assumption (A), $|\lambda_j(p^{\nu})| \geq p^{-C\nu}$ for all $1 \leq j \leq k$, for some constant $C = C(\mbf{f}) > 0$. We thus have the upper bound $|\log|\lambda_j(p^{\nu})|| \leq C\nu \log p$ for all integers $\nu \geq 2$. If $n \in \mc{N}^{(\xi)}(2x)$, in case $\nu = 1$ we have $|\log|\lambda_j(p)|| \leq \log \xi(2x)$ for any $p$ and all $1 \leq j \leq k$.  Now define 
$$
\delta_j(n) := \log|\lambda_j(n)|-\log|\tilde{\lambda}_{j,y}(n)|.
$$ 
We observe that 
\begin{align*}
\sum_{\substack{n \leq x \\ n \in \mc{N}^{(\xi)}(2x)}} |\delta_j(n)|^2 &\ll_{\mbf{f}} \sum_{n \leq x \atop n \in \mc{N}^{(\xi)}(2x)} \left(\sum_{p||n \atop p > y \text{ or } p|\Delta} |\log |\lambda_j(p)||\right)^2 + \sum_{n \leq x} \left(\sum_{\substack{p^k||n \\ k \geq 2}} k \log p\right)^2 \\
&\leq (\log \xi(2x))^2 \sum_{\substack{p,q \leq x \\ p > y \text{ or } p|\Delta \\ q > y \text{ or } q|\Delta }}\sum_{\substack{n \leq x \\ [p,q]|n}} 1 + \sum_{p^k,q^l \leq x \atop k,l \geq 2} kl (\log p)(\log q) \sum_{\substack{n \leq x \\ \nu_p(n) = k \\ \nu_q(n) = l}} 1 \\
&\ll x (\log \xi(2x))^2 \left(\sum_{p | \Delta} \frac{1}{p} + \sum_{\substack{y < p \leq x}} \frac{1}{p} \right)^2 + x (\log \xi(2x))^2 \left(\sum_{\substack{y < p \leq x}} \frac{1}{p} + \sum_{p|\Delta} \frac{1}{p}\right) \\
&+ x\left(\sum_{\substack{p \leq y \\ p^{k} \leq x \\ k \geq 2}} \frac{k \log p}{p^k}\right)^2 + x \sum_{\substack{p \leq y \\ p^k \leq x \\ k \geq 2}} \frac{k^2 (\log p)^2}{p^k} \\
&\ll_{\mbf{a},A} x(\log \log \log x)^2\left(1+\log\left(\frac{\log x}{\log y}\right)\right)^2.
\end{align*}
Since $y \geq x^{1/\log\log x}$ it follows by the union bound that
\begin{align*}
&|\{n \leq x : n \in \mc{N}^{(\xi)}(2x), |\delta_j(n)| > (\log\log\log x)^3 \text{ for some $1 \leq j \leq k$}\}|\\
& \leq (\log\log\log x)^{-6} \sum_{1 \leq j \leq k} \sum_{n \leq x \atop n \in \mc{N}^{(\xi)}(2x)} |\delta_j(n)|^2 
\ll_{\mbf{a},A} \frac{x}{(\log\log\log x)^2},
\end{align*}
and the claim follows.
\end{proof}
With the notation of Proposition \ref{prop:Sieve} we next define, for each $1 \leq j \leq k$, the quantities 
\begin{align*}
\mu_j(y) &:= \sum_{p \leq y \atop p \nmid \Delta} \frac{F_j(p)\log|\lambda_j(p)|}{p}, \ \ \ \ \ \ \ \ \sg_j(y) := \left(\sum_{p \leq y \atop p \nmid \Delta} \frac{F_j(p)(\log|\lambda_j(p)|)^2}{p}\left(1-\frac{F_j(p)}{p}\right)\right)^{1/2};
\end{align*}
Note that by Proposition \ref{prop:Sieve},
\begin{align*}
&\frac{1}{x}\sum_{\substack{n \leq x \\ n+a_j \in \mc{N}^{(\xi)}(2x) \\ \forall j}} \log|\lambda_{j,y}(n+a_j)| = \frac{1}{x}\sum_{\substack{p \leq y \\ p \nmid \Delta \\ p \in \mc{N}^{(\xi)}(2x)}} \log|\lambda_j(p)|  \sum_{\substack{n \leq x \\ p|n+a_j}} 1_{\mc{N}^{(\xi)}(2x)}(n+a_j) \\
&= \sum_{\substack{p \leq y \\ p \nmid \Delta}} \frac{F_j(p)\log|\lambda_j(p)|}{p} + O\left(\frac{1}{x}\sum_{p \leq y \atop p \in \mc{N}^{(\xi)}(2x)} |\log|\lambda_j(p)|| |r_p(x)|\right) \\
&= \mu_j(y) + O_A\left((\log \log x)^{1-A+o(1)}\right),
\end{align*}
owing to the fact that $|\lambda_j(p)| \geq (\log \log x)^{-A}$, and $y \leq x$. Thus, $\mu_j(y)$ is the mean value of $n \mapsto \log|\lambda_{j,y}(n+a_j)|$, for $n \leq x$ satisfying $n+a_j \in \mc{N}^{(\xi)}(2x)$.\\
Given non-negative integers $m_1,\ldots,m_k$, define
$$
M_{\mbf{f},\mbf{m}}(x) := \sum_{\substack{n \leq x \\ n+a_j \in \mc{N}^{(\xi)}(2x)\\ \forall j}} \prod_{1 \leq j \leq k}\left(\log|\lambda_{j,y}(n+a_j)|-\mu_j(y)\right)^{m_j}.
$$
Our goal is to prove the following result.
\begin{prop}\label{prop:Moments}
Let $x$ be large, and let $y := x^{1/\log\log\log x}$. Let $m_1,\ldots,m_k \in \mb{N} \cup \{0\}$, such that $m_j \leq \min\{\log \xi(2x),\sg_j(y)\}^{1/2}$ for all $1 \leq j \leq k$. Let $M := \max_{1 \leq j \leq k} m_j$, and assume $A > 2kM+3$. Then
\begin{align*}
M_{\mbf{f},\mbf{m}}(x) = X\prod_{1 \leq i \leq k} \frac{\Gamma(m_i+1) \sg_i(y)^{m_i}}{\Gamma(m_i/2+1)2^{m_i/2}}  \left(1_{2|m_j \forall j}+O_{k}\left((\log \xi(x))^{2}\sum_{1 \leq j\leq k} m_j^3\sg_j(y)^{-2} \right)\right).
\end{align*}
\end{prop}

\begin{proof}
We use the method of Granville and Soundararajan from \cite{GSEK}. For each $1 \leq j \leq k$ and $p \leq y$ with $p \nmid \Delta$ and $|\lambda_j(p)| \geq \xi(2x)$, we define
$$
g_{j,p}(n) := \begin{cases} \log|\lambda_j(p)|\left(1-\frac{F_j(p)}{p}\right) &\text{ if $p|n$,} \\ -\frac{F_j(p)\log|\lambda_j(p)|}{p} &\text{ if $p \nmid n$.} \end{cases}
$$
Furthermore, if $q := \prod_{p^{\nu}||q} p^{\nu}$ then set $g_{j,q}(n) := \prod_{p^{\nu}||q} g_{j,p}(n)^{\nu}$. With this notation, we observe that
\begin{align*}
M_{\mbf{f},\mbf{m}}(x) = \sum_{\substack{n \leq x \\ n+a_j \in \mc{N}^{(\xi)}(2x) \\ \forall j}} \prod_{1 \leq j \leq k} \left(\sum_{\substack{p \leq y \\ p\nmid \Delta \\ |\lambda_i(p)| \geq \xi(2x) \\ \forall i}} g_{j,p}(n+a_j)\right)^{m_j}.
\end{align*}
For convenience, for each $l \in \mb{N}$ we define 
$$
\mc{P}_{l}(y) := \{r \in \mb{N} : \Omega(r) = l, p|r \Rightarrow p \leq y, p \nmid \Delta \text{ and } |\lambda_i(p)| \geq \xi(2x) \text{ for all } 1 \leq i \leq k\}.
$$
\proofpart{1}{Employing Proposition \ref{prop:Sieve}}
\noindent Expanding the product then gives
\begin{align*}
M_{\mbf{f},\mbf{m}}(x) = \sum_{r_1 \in \mc{P}_{m_1}(y)} \cdots \sum_{r_k \in \mc{P}_{m_k}(y)} \sum_{\substack{n \leq x \\ n+a_j \in \mc{N}^{(\xi)}(2x) \\ 1 \leq j \leq k}} \prod_{1 \leq j \leq k} g_{j,r_j}(n+a_j).
\end{align*}
Write $r^{\ast}:= \text{rad}(r)$, for $r \in \mb{N}$. Let us note that if $q^{\ast} := \text{rad}(q)$ then $g_{j,q}(n) = g_{j,q}((n,q^{\ast}))$. This observation yields
$$
M_{\mbf{f},\mbf{m}}(x) = \sum_{r_1 \in \mc{P}_{m_1}(y)} \cdots \sum_{r_k \in \mc{P}_{m_k}(y)} \sum_{d_1|r_1^{\ast}}g_{1,r_1}(d_1)\cdots \sum_{d_k|r_k^{\ast}} g_{k,r_k}(d_k) \sum_{\substack{n  \leq x \\ n+a_j \in \mc{N}^{(\xi)}(2x) \\ (n+a_j,r_j^{\ast}) = d_j \\ \forall j}} 1.
$$
Note that $(r_i,r_j) = 1$ for all $i \neq j$ immediately, given that any $p|(r_i,r_j)$ would divide $\Delta$, while $(r_j,\Delta) = 1$ for all $j$ by construction.
By M\"{o}bius inversion, the inner sum is
$$
\sum_{\substack{n  \leq x \\ n+a_j \in \mc{N}^{(\xi)}(2x) \\ d_j|(n+a_j)}} \prod_{1 \leq j \leq k} \sum_{e_j|r_j^{\ast}/d_j} \mu(e_j)1_{e_j|(n+a_j)/d_j} = \sum_{e_1d_1|r_1^{\ast}} \mu(e_1) \cdots \sum_{e_kd_k|r_k^{\ast}} \mu(e_k) \sum_{\substack{n  \leq x \\ n+a_j \in \mc{N}^{(\xi)}(2x) \\ e_jd_j|(n+a_j)}} 1,
$$
so that
$$
M_{\mbf{f},\mbf{m}}(x) = \mathop{\sum_{r_1 \in \mc{P}_{m_1}(y)} \cdots \sum_{r_k \in \mc{P}_{m_k}(y)}}_{(r_i,r_j) = 1 \ \forall i \neq j} \sum_{e_1d_1|r_1^{\ast}}\mu(e_1)g_{1,r_1}(d_1)\cdots \sum_{e_kd_k|r_k^{\ast}} \mu(e_k)g_{k,r_k}(d_k) \sum_{\substack{n  \leq x \\ n+a_j \in \mc{N}^{(\xi)}(2x) \\ e_jd_j|(n+a_j)}} 1.
$$
Since $(r_i,r_j) = 1$ for all $i \neq j$ we may apply Proposition \ref{prop:Sieve} (noting that $y = x^{o(1)}$) to give
\begin{align*}
M_{\mbf{f},\mbf{m}}(x) = \mathop{\sum_{r_1 \in \mc{P}_{m_1}(y)} \cdots \sum_{r_k \in \mc{P}_{m_k}(y)}}_{(r_i,r_j) = 1 \forall i \neq j} &\sum_{e_1d_1|r_1^{\ast}} \mu(e_1)g_{1,r_1}(d_1)\cdots \sum_{e_kd_k|r_k^{\ast}} \mu(e_k)g_{k,r_k}(d_k) \\
&\cdot \frac{X}{d_1e_1\cdots d_ke_k} \left(F_1(d_1e_1)\cdots F_k(d_ke_k) + O_A( (\log\log x)^{1-A})\right) \\
&=: \tilde{\mc{M}} + \tilde{\mc{E}},
\end{align*}
where $X$ is as in the statement of Proposition \ref{prop:Sieve}. \\
We first treat the main term $\tilde{\mc{M}}$. We obtain
\begin{align*}
\tilde{\mc{M}} = X\mathop{\sum_{r_1 \in \mc{P}_{m_1}(y)} \cdots \sum_{r_k \in \mc{P}_{m_k}(y)}}_{\substack{(r_i,r_j) = 1 \\ \forall i \neq j}} \prod_{1 \leq j \leq k} \sum_{e_jd_j|r_j^{\ast}} \frac{\mu(e_j)F_j(e_j)}{e_j} \cdot \frac{g_{j,r_j}(d_j)F_j(d_j)}{d_j}.
\end{align*}
The product over $j$ is
\begin{align*}
&\prod_{1 \leq j \leq k} \sum_{b_je_jd_j = r_j^{\ast}} \frac{\mu(e_j)F_j(e_j)}{e_j}\prod_{p|b_je_j \atop p^{\nu}||r_j} \left(-\frac{F_j(p)\log|\lambda_j(p)|}{p}\right)^{\nu} \cdot \frac{F_j(d_j)}{d_j}\prod_{p|d_j \atop p^{\nu}||r_j} (\log|\lambda_j(p)|)^{\nu} \left(1-\frac{F_j(p)}{p}\right)^{\nu} \\
&= \prod_{1 \leq j \leq k} \prod_{p|r_j^{\ast} \atop p^{\nu}||r_j} \left(\left(1-\frac{F_j(p)}{p}\right)\left(-\frac{F_j(p)\log|\lambda_j(p)|}{p}\right)^{\nu} + (\log|\lambda_j(p)|)^{\nu}\frac{F_j(p)}{p}\left(1-\frac{F_j(p)}{p}\right)^{\nu} \right) \\
&=: \prod_{1 \leq j \leq k} G_j(r_j).
\end{align*}
Note that for each $1 \leq j \leq k$, any $p \leq y$ with $p\nmid \Delta$ and $\lambda_j(p) \neq 0$, and $\alpha \geq 1$,
\begin{align} \label{eq:Gjexp}
G_j(p^{\alpha}) &= (\log|\lambda_j(p)|)^{\alpha} \left(\frac{F_j(p)}{p} \left(1-\frac{F_j(p)}{p}\right)^{\alpha} + \left(1-\frac{F_j(p)}{p}\right)\left(-\frac{F_j(p)}{p}\right)^{\alpha}\right) \nonumber \\
&= (\log|\lambda_j(p)|)^{\alpha} \frac{F_j(p)}{p}\left(1-\frac{F_j(p)}{p}\right) \left(\left(1-\frac{F_j(p)}{p}\right)^{\alpha-1} - \left(-\frac{F_j(p)}{p}\right)^{\alpha-1}\right),
\end{align}
so that $\prod_{1 \leq j \leq k} G_j(r_j) = 0$ unless $r_j$ is square-full for all $1 \leq j \leq k$, i.e., $\nu_p(r_j) \geq 2$ for all $p|r_j$ and $1 \leq j \leq k$. We thus obtain
$$
\tilde{\mc{M}} = X\mathop{\sum_{r_1 \in \mc{P}_{m_1}(y) \atop r_1 \text{ square-full}} \cdots \sum_{r_k \in \mc{P}_{m_k}(y) \atop r_k \text{ square-full}}}_{(r_i,r_j) = 1 \forall i \neq j} \prod_{1 \leq j \leq k} G_j(r_j).
$$
\proofpart{2}{A variant of $\tilde{\mc{M}}$}
\noindent Writing $\mc{M}$ to denote the same expression \emph{without} the constraint $(r_i,r_j) = 1$, we obtain
\begin{align*}
\mc{M} &= X\prod_{1 \leq j \leq k}\sum_{r_j \in \mc{P}_{m_j}(y) \atop r_j \text{ square-full}} G_j(r_j) = \prod_{1 \leq j \leq k} \sum_{1 \leq l_j \leq m_j/2} \sum_{\substack{q_1 < \cdots <  q_{l_j} \leq y \\ q_i \text{ distinct }, q_i \nmid \Delta \\ |\lambda_j(q_i)| \geq \xi(2x) \\ \forall i,j}} \sum_{\alpha_1,\ldots,\alpha_{l_j}\geq 2 \atop \alpha_1 + \cdots + \alpha_{l_j} = m_j} \frac{m_j!}{l_j!\alpha_1!\cdots \alpha_{l_j}!} \prod_{1 \leq i \leq l_j} G_j\left(q_i^{\alpha_i}\right) \\
&=: \prod_{1 \leq j \leq k} \sum_{1 \leq l_j \leq m_j/2} T_{j,l_j}(x;y).
\end{align*}
Fix $1 \leq j \leq k$ and $1 \leq l_j < m_j/2$ for the moment. Note from \eqref{eq:Gjexp} and the obvious bound 
$$
|(1-u)^r -(-u)^r| \leq (1-u)^r + u^r \leq 1-u + u = 1,
$$
for any $u \in [0,1]$ and $r \geq 1$, that (after dispensing with the distinctness condition)
\begin{align*}
&T_{j,l_j}(x;y) \leq \frac{m_j!}{l_j!}\sum_{\alpha_1 + \cdots + \alpha_{l_j} = m_j \atop \alpha_j \geq 2} \prod_{1 \leq i \leq l_j} \frac{1}{\alpha_j!} \left(\sum_{\substack{q_i \leq y \\ q_i \nmid \Delta \\ |\lambda_j(q_i)| \geq \xi(2x) \forall i,j}} \frac{|\log|\lambda_j(q_i)||^{\alpha_j}F_j(q_i)}{q_i}\left(1-\frac{F_j(q_i)}{q_i}\right)\right) \\
&\leq \frac{m_j!}{l_j!2^{l_j}}(\log \xi(2x))^{m_j-2l_j}\sg_j(y)^{2l_j}  |\{(\alpha_1,\ldots,\alpha_{l_j}) \in \mb{N}^{l_j} : \alpha_1+\cdots + \alpha_{l_j} = m_j, \alpha_1,\ldots, \alpha_{l_j} \geq 2\}| \\
&= \frac{m_j!}{l_j!2^{l_j}}\binom{m_j-l_j-1}{l_j-1} (\log \xi(2x))^{m_j-2l_j} \sg_j(y)^{2l_j}.
\end{align*}
Comparing these bounds for $l_j$ and $l_j+1$, we see that when $m_j \leq (\sg_j(y)^2/\log \xi(2x))^{1/3}$ we obtain
$$
\max_{1 \leq l_j < m_j/2} T_{j,l_j}(x;y) \ll \frac{m_j^3(\log \xi(2x))^2}{\sg_j(y)^2} \cdot \frac{\Gamma(m_j+1)}{2^{m_j/2}\Gamma(m_j/2+1)}\sg_j(y)^{m_j}.
$$
Thus, provided that $2|m_j$ for all $j$, we relegate to the error term all but the arrangements with $l_j = m_j/2$, and $\alpha_i = 2$ for all $1 \leq i \leq m_j/2$, the number of which is precisely $\frac{m_j!}{2^{m_j/2}(m_j/2)!}$. It follows, therefore, that 
\begin{align*}
\mc{M} &= X\prod_{1 \leq j \leq k} 1_{2|m_j} \frac{m_j!}{2^{m_j/2}(m_j/2)!}  \sum_{\substack{q_1 < \ldots < q_{m_j/2} \leq y \\ q_i \nmid \Delta \\ |\lambda_j(q_i)| \geq \xi(2x) \forall i,j}} \left(\prod_{1 \leq i \leq m_j/2} G_j(q_i^2)\right) \\
&+ O_k\left(x(\log \xi(2x))^2 \left(\prod_{1 \leq j \leq k} \frac{\Gamma(m_j+1)}{2^{m_j/2}\Gamma(m_j/2+1)}\sg_j(y)^{m_j} \right) \sum_{1 \leq j' \leq k} m_{j'}^3\sg_{j'}(y)^{-2}\right).
\end{align*}
Removing the distinctness condition for $q_i$ in the main term of the above estimate when some $m_j > 2$ incurs (up to relabeling) an error of size
\begin{align*}
&\ll_k x\prod_{1 \leq j \leq k} \frac{\Gamma(m_j+1)}{2^{m_j/2}\Gamma(m_j/2+1)}\left|\sum_{\substack{q_1,\ldots,q_{m_j/2} \leq y \\ q_{m_j/2} = q_{m_j/2-1} \\ |\lambda_j(q_i)| \geq \xi(2x) \forall i,j}} G_j(q_{m_j/2}^2)^2\prod_{1 \leq i \leq m_j/2-2} G_j(q_i^2)\right| \\
&\ll x(\log \xi(2x))^4 \left(\prod_{1 \leq j \leq k} \frac{\Gamma(m_j+1)}{2^{m_j/2}\Gamma(m_j/2+1)}\sg_j(y)^{m_j}\right) \sum_{1 \leq j' \leq k} \sg_{j'}(y)^{-4}.
\end{align*}
We thus conclude that
\begin{align}
\mc{M} &= X\prod_{1 \leq j \leq k} 1_{2|m_j}\frac{m_j!}{2^{m_j/2}(m_j/2)!}\left(\sum_{\substack{q_i \leq y \\ q_i \nmid \Delta \\ |\lambda_j(q_i)| \geq \log \xi(2x) \forall i,j}} (\log|\lambda_j(q_i)|)^2 \frac{F_j(q_i)}{q_i}\left(1-\frac{F_j(q_i)}{q_i}\right) \right)^{m_j/2} \nonumber \\
&+ O_k\left((\log \xi(2x))^{2}\left(\prod_{1 \leq j \leq k} \frac{\Gamma(m_j+1)}{2^{m_j/2}\Gamma(m_j/2+1)}\sg_j(y)^{m_j}\right)\sum_{1 \leq j' \leq k} m_{j'}^3\sg_{j'}(y)^{-2}\right) \nonumber \\
&= X\prod_{1 \leq j \leq k} \frac{\Gamma(m_j+1)}{2^{m_j/2}\Gamma(m_j/2+1)} \sg_j(y)^{m_j} \left(1_{2|m_j \forall j}+O_k\left((\log \xi(2x))^2 \sum_{1 \leq j' \leq k} m_{j'}^3\sg_{j'}(y)^{-2}\right)\right). \label{eq:momBadErr}
\end{align}
\proofpart{3}{Comparing $\mc{M}$ with $\tilde{\mc{M}}$} 
\noindent Note that
$$
|\mc{M}-\tilde{\mc{M}}| \ll x\mathop{\sum_{r_1 \in \mc{P}_{m_1}(y) \atop r_1 \text{ square-full}} \cdots \sum_{r_k \in \mc{P}_{m_k}(y) \atop r_k \text{ square-full}}}_{\exists i \neq j: (r_i,r_j) > 1} \prod_{1 \leq j \leq k} |G_j(r_j)|.
$$
Fix $1 \leq i < j \leq k$, such that $(r_i,r_j) > 1$ (assuming such a pair exists). Then there is some square-full $d$ such that $d|r_{i}$ and $d|r_j$, with $\Omega(d) \leq \min\{m_i,m_j\}$, and each of its prime factors $p$ is $\leq y$ and satisfies $|\lambda_j(p)| \geq \xi(2x)$ for all $1 \leq j \leq k$. It follows by the union bound and an earlier calculation that
\begin{align*}
&|\mc{M}-\tilde{\mc{M}}| \ll x\sum_{1 \leq i < j \leq k} \left(\prod_{1 \leq j' \leq k \atop j' \neq i,j} \sum_{r_{j'} \in \mc{P}_{m_{j'}}(y) \atop r_{j'} \text{ square-full}}|G_{j'}(r_{j'})|\right) \sum_{\substack{P^+(d) \leq y \\ 1 \leq \Omega(d) \leq \min\{m_i,m_j\} \\ d \text{ square-full}}} \sum_{r_i \in \mc{P}_{m_i}(y) \atop d|r_i} |G_i(r_i)| \sum_{r_j \in \mc{P}_{m_j}(y) \atop d|r_j} |G_j(r_j)| \\
&\ll x\sum_{1 \leq i < j \leq k} \left(\prod_{1 \leq j' \leq k \atop j' \neq i,j} \frac{\Gamma(m_{j'}+1)\sg_{j'}(y)^{m_{j'}}}{2^{m_{j'}/2}\Gamma(m_{j'}/2+1)}\right)\sum_{2 \leq r \leq \min\{m_i,m_j\}} \sum_{1 \leq s \leq r/2} \sum_{\alpha_1+\cdots + \alpha_s = r \atop \alpha_i \geq 2} \frac{r!}{s!\alpha_1!\cdots \alpha_s!} \\
&\cdot \sum_{q_1 < \cdots < q_s \leq y} \sum_{r_i' \in \mc{P}_{m_i-r}(y) \atop r_i' \prod_t q_t^{\alpha_t} \text{ square-full}} \left|G_i\left(r_i'\prod_{1\leq t\leq s} q_t^{\alpha_t}\right)\right| \sum_{r_j' \in \mc{P}_{m_j-r}(y) \atop r_j' \prod_t q_t^{\alpha_t} \text{ square-full}} \left|G_j\left(r_j'\prod_{1\leq t\leq s} q_t^{\alpha_t}\right)\right|.
\end{align*}
Fixing $q_1 < \ldots < q_s$ and writing $r_i' = r_i''\prod_{1 \leq t \leq s} q_t^{\gamma_t}$, where $(r_i'',\prod_t q_t) = 1$, $\gamma_t \geq 0$ for all $1 \leq t \leq s$ and $\gamma_1 + \cdots + \gamma_s \leq m_i-r$, the sum over $r_i'$ can be bounded by
\begin{align*}
&\sum_{\gamma = \gamma_1 + \cdots + \gamma_s \leq m_i-r \atop \gamma_t \geq 0} \prod_{1 \leq t \leq s} |G_i(q_t^{\alpha_t + \gamma_t})| \sum_{r_i'' \in \mc{P}_{i,m_i-r-\gamma}(y) \atop r_i'' \text{ square-full}} |G_i(r_i'')| \\
&\ll \sum_{\gamma = \gamma_1 + \cdots + \gamma_s \leq m_i-r \atop \gamma_t \geq 0} \prod_{1 \leq t \leq s} |G_i(q_t^{\alpha_t + \gamma_t})| \frac{\Gamma(m_i-r-\gamma+1)}{2^{\frac{1}{2}(m_i-r-\gamma)}\Gamma((m_i-r-\gamma)/2+1)} \sg_i^{m_i-r-\gamma}(y).
\end{align*}
Given $\gamma$, Stirling's formula implies that
\begin{align*}
\frac{\Gamma(m_i-r-\gamma+1)}{2^{(m_i-r-\gamma)/2}\Gamma((m_i-r-\gamma)/2+1)} &\ll \left(\frac{e}{2m_i}\right)^{(r+\gamma)/2}\left(1-\frac{r+\gamma}{m_j}\right)^{m_i/2} \frac{\Gamma(m_i+1)}{2^{m_i/2}\Gamma(m_i/2+1)} \\
&\ll m_i^{-(r+\gamma)/2}\frac{\Gamma(m_i+1)}{2^{m_i/2}\Gamma(m_i/2+1)}.
\end{align*}
Applying this with the sum over $r_i'$ and a similar bound for the sum in $r_j'$, then inserting both above, we obtain
\begin{align*}
&|\mc{M}-\tilde{\mc{M}}|\ll x\left(\prod_{1 \leq j' \leq k} \frac{\Gamma(m_{j'}+1) \sg_{j'}(y)^{m_{j'}}}{2^{m_{j'}/2}\Gamma(m_{j'}/2+1)}\right)\sum_{1 \leq i < j \leq k} \sum_{2 \leq r \leq \min\{m_i,m_j\}} \frac{r!}{(\sg_i(y)\sg_j(y))^r} \sum_{1 \leq s \leq r/2} \frac{1}{s!} \\
&\cdot \sum_{\alpha_1+\cdots + \alpha_s = r \atop \alpha_i \geq 2} \frac{1}{\alpha_1!\cdots \alpha_s!} \sum_{\gamma \leq m_i-r} \frac{1}{m_i^{(r+\gamma)/2}\sg_i(y)^{\gamma}}\sum_{\delta \leq m_j-r} \frac{1}{m_j^{(r+\delta)/2}\sg_j(y)^{\delta}} \sum_{\gamma_1 + \cdots + \gamma_s = \gamma \atop \gamma_t \geq 0} \sum_{\delta_1 + \cdots + \delta_s = \delta \atop \delta_t \geq 0} \\
&\cdot \prod_{1 \leq t \leq s} \left(\sum_{\substack{q_t \leq y \\ \lambda_i(q_t)\lambda_j(q_t) \neq 0}} |G_i(q_t^{\alpha_t + \gamma_t})||G_j(q_t^{\alpha_t + \delta_t})|\right).
\end{align*}
The product over $1 \leq t \leq s$ is
$$
\leq \prod_{1 \leq t \leq s} \sum_{q_t \leq y \atop |\lambda_i(q_t)||\lambda_j(q_t)| \neq 0} \frac{|\log|\lambda_i(q_t)||^{\alpha_t+\gamma_t}|\log|\lambda_j(q_t)||^{\alpha_t+\delta_t}}{q_t^2} \leq C^s (\log \xi(2x))^{2r+\gamma+\delta}
$$
for some constant $C > 0$. 
For $a \in \mb{N}$, set 
$$
p_s(a) := |\{(\beta_1,\ldots,\beta_s) \in (\mb{N} \cup \{0\})^s : \beta_1 + \cdots + \beta_s = a\}| \ll \frac{a^{s-1}}{(s-1)!}.
$$ 
As $s \leq r/2$ and $\gamma < m_i$ we have $p_s(\gamma)m_i^{-r/2} \ll 1$, and similarly $p_s(\delta)m_j^{-r/2} \ll 1$. As such, for each $1 \leq s \leq r/2$,
\begin{align*}
&\sum_{\gamma \leq m_i-r} \frac{p_s(\gamma) (\log \xi(2x))^{\gamma}}{m_i^{(r+\gamma)/2}\sg_i(y)^{\gamma}}\sum_{\delta \leq m_j-r} \frac{p_s(\delta)(\log \xi(2x))^{\delta}}{m_j^{(r+\delta)/2}\sg_j(y)^{\delta}} \ll \sum_{\gamma \leq m_i-r} \left(\frac{\log \xi(2x)}{m_i^{1/2}\sg_i(y)}\right)^{\gamma} \sum_{\delta \leq m_j-r} \left(\frac{\log \xi(2x)}{m_j^{1/2}\sg_j(y)}\right)^{\delta} \ll 1.
\end{align*}
We thus find that
\begin{align*}
|\mc{M}-\tilde{\mc{M}}| &\ll x\left(\prod_{1 \leq j' \leq k} \frac{\Gamma(m_{j'}+1) \sg_{j'}(y)^{m_{j'}}}{2^{m_{j'}/2}\Gamma(m_{j'}/2+1)}\right)\sum_{1 \leq i < j \leq k} \sum_{2 \leq r \leq \min\{m_i,m_j\}} \frac{r!(\log \xi(2x))^{2r}}{(\sg_i(y)\sg_j(y))^r} \sum_{1 \leq s \leq r/2} \frac{C^s}{s!} \\
&\cdot \sum_{\alpha_1+\cdots + \alpha_s = r \atop \alpha_i \geq 2} \frac{1}{\alpha_1!\cdots \alpha_s!} \sum_{\gamma \leq m_i-r} \frac{p_s(\gamma)(\log \xi(2x))^{\gamma}}{m_i^{(r+\gamma)/2}\sg_i(y)^{\gamma}}\sum_{\delta \leq m_j-r} \frac{p_s(\delta)(\log \xi(2x))^{\delta}}{m_j^{(r+\delta)/2}\sg_j(y)^{\delta}}\\
&\ll x\left(\prod_{1 \leq j' \leq k} \frac{\Gamma(m_{j'}+1) \sg_{j'}(y)^{m_{j'}}}{2^{m_{j'}/2}\Gamma(m_{j'}/2+1)}\right)\sum_{1 \leq i < j \leq k} \sum_{2 \leq r \leq \min\{m_i,m_j\}} \left(\frac{(\log \xi(2x))^2}{\sg_i(y)\sg_j(y)}\right)^r \\
&\cdot \sum_{1 \leq s \leq r/2} \frac{C^s}{s!} \sum_{\alpha_1+\cdots + \alpha_s = r \atop \alpha_i \geq 2} \frac{r!}{\alpha_1!\cdots \alpha_s!} \\
&\leq e^C x\left(\prod_{1 \leq j' \leq k} \frac{\Gamma(m_{j'}+1) \sg_{j'}(y)^{m_{j'}}}{2^{m_{j'}/2}\Gamma(m_{j'}/2+1)}\right)\sum_{1 \leq i < j \leq k} \sum_{2 \leq r \leq \min\{m_i,m_j\}} \frac{(m_im_j)^{r/2}(\log \xi(2x))^{2r}}{(\sg_i(y)\sg_j(y))^r},
\end{align*}
where we used the fact that the terms $r!/(\alpha_1!\cdots \alpha_s!)$ are multinomial coefficients, whose sum is $\leq s^r \leq r^r$, and $r \leq \min\{m_i,m_j\} \leq (m_im_j)^{1/2}$. Since $m_l \leq \sg_l(y)^{1/2}$ for all $1 \leq l \leq k$ we finally obtain, after removing the condition $i \neq j$,
$$
|\mc{M}-\tilde{\mc{M}}| \ll x(\log \xi(2x))^4\left(\prod_{1 \leq j \leq k} \frac{\Gamma(m_j+1)}{2^{m_j/2}\Gamma(m_j/2+1)} \sg_j(y)^{m_j} \right) \cdot \left(\sum_{1 \leq l \leq k} m_l\sg_l(y)^{-2}\right)^2.
$$
\proofpart{4}{Treatment of $\tilde{\mc{E}}$ and Conclusion} 
\noindent We now estimate $\tilde{\mc{E}}$. This can be bounded by
\begin{align*}
|\tilde{\mc{E}}| &\ll \frac{x}{(\log \log x)^{A-1}} \prod_{1 \leq j \leq k} \sum_{r_j \in \mc{P}_{j,m_j}(y)} \sum_{e_jd_j|r_j^{\ast}} \frac{|g_{j,r_j}(d_j)|}{e_jd_j} \\
&\ll \frac{x(\log \xi(2x))^{kM}}{(\log \log x)^{A-1}} \prod_{1 \leq j \leq k} \sg_j(y)^{2m_j} \ll \frac{x}{(\log \log x)^{A-kM -1 + o_A(1)}} \prod_{1 \leq j \leq k} \sg_j(y)^{m_j}
\end{align*}
given that $\xi(x) \geq (\log \log x)^{-A}$ and $M \ll_A (\log \log \log x)^{1/2}$. Thus, as $A > 2kM+3$ we obtain
\begin{align*}
M_{\mbf{f},\mbf{m}}(x) &= X \prod_{1 \leq j \leq k} \frac{\Gamma(m_j+1)}{2^{m_j/2}\Gamma(m_j/2+1)}\sg_j(y)^{m_j} \left(1_{2|m_j \forall j}+O_k\left((\log \xi(2x))^{4}\sum_{1 \leq j \leq k} m_j^3\sg_j(y)^{-2}\right)\right) \\
&+ O_A\left(\frac{x}{(\log \log x)^{A/2}}\right),
\end{align*}
which completes the proof.
\end{proof}

\begin{proof}[Proof of Proposition \ref{prop:EKCusp}]
Let $x$ be large, and let $y := x^{\frac{1}{\log\log\log x}}$. Let $A > 1$ be large and set $\xi(x) := (\log\log x)^{-A}$ as above. Put $c := \frac{1}{2}(1+\pi^2/6)$. For each $1 \leq j \leq k$, Lemma 11 of \cite{LRS} yields 
\begin{align*}
\mu_j(y) &= \sum_{\substack{p \leq y \\ p \nmid \Delta \\ |\lambda_i(p)| \geq \xi(2x) \forall i}} \frac{\log|\lambda_j(p)|F_j(p)}{p} = -\frac{1}{2}\log\log x + O_A(\log\log \log x) \\
\sg_j(y)^2 &= \sum_{\substack{p \leq y \\ p \nmid \Delta \\ |\lambda_i(p)| \geq \xi(2x) \forall i}} \frac{(\log |\lambda_j(p)|)^2F_j(p)}{p}\left(1-\frac{F_j(p)}{p}\right) = c \log\log x + O_A(\log\log\log x).
\end{align*}
For each $n \leq x$ with $n \in \mc{N}$ and $1 \leq j \leq k$ define
\begin{align*}
X_j(n) &:= \frac{\log|\lambda_j(n)| + \frac{1}{2}\log\log n}{\sqrt{c \log\log n}},
\end{align*}
and for $n \in \mc{N}^{(\xi)}(2x)$ define
\begin{align*}
X_{j,y}(n) &:= \frac{\log|\tilde{\lambda}_{j,y}(n)| - \mu_j(y)}{\sg_j(y)}.
\end{align*}
Let $I_1,\ldots, I_k$ be closed non-empty intervals. Let us define
$$
N_{\mbf{I},\mbf{a}}(x) := |\{n \leq x : n+a_j \in \mc{N}, X_j(n+a_j) \in I_j,1 \leq j \leq k\}|
$$
as well as
$$
N_{\mbf{I},\mbf{a}}^{(\xi)}(x;y) := |\{n \leq x : n+a_j \in \mc{N}^{(\xi)}(2x), X_{j,y}(n+a_j) \in I_j, 1 \leq j \leq k\}|.
$$
By Lemma \ref{lem:approx}, we have
\begin{align} \label{eq:XiCond}
N_{\mbf{I},\mbf{a}}(x) &= |\{n \leq x : n+a_j \in \mc{N}^{(\xi)}(2x), X_j(n+a_j) \in I_j \forall 1 \leq j \leq k\}| + O\left(\frac{x}{(\log\log x)^{A-1}}\right).
\end{align}
It follows from Lemma \ref{lem:assump} that for all but $o(x)$ integers $\sqrt{x} < n \leq x$ with $n \in \mc{N}^{(\xi)}(2x)$ we have
$$
\max_{1 \leq j \leq k} |X_j(n)-X_{j,y}(n)| \ll_A \frac{(\log\log\log x)^3}{\sqrt{\log\log x}}.
$$
Thus, in particular, if $x$ is large enough then we can find intervals $I_j^- \subseteq I_j \subseteq I_J^+$, such that $\text{meas}(I_j^+\bk I_j^-) \ll \frac{(\log\log\log x)^3}{\sqrt{\log\log x}}$ and
\begin{align}\label{eq:ULBds}
N_{\mbf{I}^-,\mbf{a}}^{(\xi)}(x;y) &\leq |\{n \leq x : n+a_j \in \mc{N}^{(\xi)}(2x), X_j(n+a_j) \in I_j \forall 1 \leq j \leq k\}| + o(x) \nonumber \\
&\leq N_{\mbf{I}^+,\mbf{a}}^{(\xi)}(x;y) + o(x).
\end{align}
We claim now that for any collection of non-empty closed intervals $J_1,\ldots,J_k$ we have
\begin{equation}\label{eq:claim}
\lim_{x \ra \infty} \frac{1}{X} N_{\mbf{J},\mbf{a}}^{(\xi)}(x;y) = \frac{1}{(2\pi)^k} \int_{\mc{B}(\mbf{J})} e^{-\frac{1}{2}\mbf{x} \cdot \mbf{x}} d\mbf{x},
\end{equation}
where $\mbf{x} \cdot \mbf{x} = \sum_{1 \leq j \leq k} x_j^2$ denotes the usual scalar product, and $X = X(x)$ is as in Proposition \ref{prop:Sieve}. Assuming this, it will follow from combining \eqref{eq:XiCond} and \eqref{eq:ULBds} that
$$
\frac{1}{(2\pi)^k} \int_{\mc{B}(\mbf{I}^-)} e^{-\frac{1}{2}\mbf{x} \cdot \mbf{x}} d\mbf{x} \leq \frac{1}{X}N_{\mbf{I},\mbf{a}}(x) + o_{\mbf{a},k}(1) \leq \frac{1}{(2\pi)^k} \int_{\mc{B}(\mbf{I}^+)} e^{-\frac{1}{2}\mbf{x} \cdot \mbf{x}} d\mbf{x} + o_{\mbf{a},k}(1),
$$
and given that $\text{meas}(I_j^+\bk I_j^-) = o(1)$ and the absolute continuity of the Gaussian distribution,
$$
\frac{1}{X}N_{\mbf{I},\mbf{a}}(x) = \frac{1}{(2\pi)^k} \int_{\mc{B}(\mbf{I})} e^{-\frac{1}{2}\mbf{x} \cdot \mbf{x}} d\mbf{x} +o_{\mbf{a},k}(1),
$$
as required.\\
It remains to check \eqref{eq:claim}. Let $Z_1,\ldots, Z_k$ be pairwise independent Gaussian random variables with mean 0 and variance 1. By the method of moments, it is sufficient to show that for any non-negative integers $m_1,\ldots,m_k$, 
\begin{equation}\label{eq:MomCheck}
\frac{1}{X}\sum_{n \leq x \atop n+a_j \in \mc{N}_j^{(\xi)}(x) \forall j} X_{1,y}(n+a_1)^{m_1} \cdots X_{k,y}(n+a_k)^{m_k} = \mbf{E}[Z_1^{m_1}\cdots Z_k^{m_k}] + o(1),
\end{equation}
as $x \ra \infty$. By independence, it is well known that
$$
\mbf{E}[Z_1^{m_1}\cdots Z_k^{m_k}] = \prod_{1 \leq j \leq k} \mbf{E}[Z_j^{m_j}] = \prod_{1 \leq j \leq k} 1_{2|m_j}\frac{\Gamma(m_j+1)}{2^{m_j/2}\Gamma(m_j/2+1)},
$$
Proposition \ref{prop:Moments} implies, on the other hand, that
\begin{align*}
&\frac{1}{X} \sum_{n \leq x \atop n+a_j \in \mc{N}_j^{(\xi)}(2x) \forall j} X_{1,y}(n+a_1)^{m_1} \cdots X_{k,y}(n+a_k)^{m_k} \\
&= \prod_{1 \leq j \leq k} \frac{\Gamma(m_j+1)}{2^{m_j/2}\Gamma(m_j/2+1)}\left(1_{2|m_i \forall i}  + O_{k,A}\left(M^3\frac{(\log \log\log x)^4}{\log \log x}\right)\right),
\end{align*}
where $M := \max_{1 \leq j \leq k} m_j$, provided that $A > 2kM+3$. Since $A$ was arbitrarily large, as $x \ra \infty$ we see therefore that \eqref{eq:MomCheck} holds, and the proof of the proposition is complete.
\end{proof}

\section{Proof of Proposition~\ref{prop:EKwithSigns}} \label{sec:signs}
\subsection{Fourier analytic reductions} As in the proof of Proposition \ref{prop:EKCusp}, we can find tuples of intervals $\mbf{I}^{(1)},\mbf{I}^{(2)}$ such that
\begin{align} \label{eq:smoothSign}
N_{\mbf{a}}(x;\mbf{\eps},\mbf{I}) &\geq |\{n \leq x : n+a_j \in \mc{N}^{(\xi)}(2x), \sg_j(n+a_j) = \eps_j, X_{j,y}(n+a_j) \in I_j^{(1)}, 1 \leq j \leq k\}| + o(x) \nonumber \\
&=: N_{\mbf{a},y}(x;\mbf{\eps},\mbf{I}^{(1)}) + o(x),
\end{align}
and similarly $N_{\mbf{a}}(x;\mbf{\eps},\mbf{I}) \leq N_{\mbf{a},y}(x;\mbf{\eps},\mbf{I}^{(2)}) + o(x)$, where $y = x^{1/\log\log\log x}$ and $\xi(x) = (\log\log x)^{-A}$ for some $A > 1$. 
Let $\delta=\delta(A,k,\mbf{a},\mbf{I})$ be a parameter to be chosen, and let $\tilde{\mbf{I}} \in \{\mbf{I}^{(1)},\mbf{I}^{(2)}\}$. For each $1 \leq j \leq k$ we select a weight function $w_j \in \mc{C}^{\infty}_c(\mb{R})$ (depending on $\delta$) that satisfies: 
\begin{enumerate}[label = (\roman*)]
\item $\text{supp}(w_j) \subseteq \tilde{I}_{j,\delta} := \{x \in \mb{R} : (x-\delta,x+\delta) \cap \tilde{I}_j \neq \emptyset\}$,  \\
\item $0 \leq w_j \leq 1$ with $w_j \equiv 1$ on $\tilde{I}_j$, and \\
\item for all $B > 0$,
\begin{equation} \label{eq:FTbd}
|\hat{w}_j(\xi)| \ll_{\mbf{I},B} (1+|\xi|^2)^{-B}.
\end{equation}
\end{enumerate}
We observe that, as a consequence of Proposition \ref{prop:EKCusp},
\begin{align*}
&N_{\mbf{a},y}(x;\mbf{\eps},\tilde{\mbf{I}}) = \sum_{n \leq x} \prod_{1 \leq j \leq k} 1_{\mc{N}^{(\xi)}(2x)}(n+a_j) 1_{\sg_j(n+a_j) = \e_j} w_j(X_{j,y}(n+a_j)) \\
&+ O\left(\sum_{1 \leq j \leq k} \sum_{n \leq x} \left(\prod_{1 \leq j' \leq k \atop j' \neq  j} 1_{\mc{N}^{(\xi)}(2x)}(n+a_{j'}) 1_{\tilde{I}_{j',\delta}}(X_{j,y}(n+a_{j'}))\right) 1_{\mc{N}^{(\xi)}(2x)}(n+a_j) 1_{\tilde{I}_{j,\delta} \bk \tilde{I}_j}(X_{j,y}(n+a_j))\right) \\
&= \sum_{n \leq x} \prod_{1 \leq j \leq k} 1_{\mc{N}^{(\xi)}(2x)}(n+a_j) 1_{\sg_j(n+a_j) = \e_j} w_j(X_{j,y}(n+a_j)) + O_k(\delta x) \\
&=: \tilde{N}_{\mbf{a},y}(x;\mbf{\eps},\tilde{I}) + O_k(\delta x).
\end{align*}
Using the identity 
\begin{equation}\label{eq:pmIden}
1_{\sg_j(n) = \e_j} = \frac{1}{2}(1+\e_j \sg_j(n))
\end{equation}
for each $n \in \mb{N}$ and $1 \leq j \leq k$, and Fourier inversion, we obtain
\begin{align*}
&\tilde{N}_{\mbf{a},y}(x;\mbf{\eps},\tilde{\mbf{I}}) = \sum_{n \leq x} \prod_{1 \leq j \leq k} 1_{\mc{N}^{(\xi)}(2x)}(n+a_j) 1_{\sg_j(n+a_j) = \e_j} w_j(X_{j,y}(n+a_j)) \\
&= 2^{-k} \sum_{\mbf{\nu} \in \{0,1\}^k} \left(\prod_{1 \leq i \leq k} \e_i^{\nu_i}\right) \int_{\mb{R}^k} W(\mbf{u}) \left(\sum_{n \leq x} \prod_{1 \leq j \leq k} 1_{\mc{N}^{(\xi)}(2x)}(n+a_j) \sg_j^{\nu_j}(n+a_j) e(u_j X_{j,y}(n+a_j))\right) d\mbf{u},
\end{align*}
where we defined
$$
W(\mbf{u}) := \prod_{1 \leq j \leq k} \hat{w}_j(u_j).
$$
Let $1 \leq Z\leq \log\log \log x$ be a parameter of our choosing. By the union bound, we find that
\begin{align*}
&\int_{\|\mbf{u}\|_{\infty} > Z} W(\mbf{u}) \left(\sum_{n \leq x} \prod_{1 \leq j \leq k} 1_{\mc{N}^{(\xi)}(2x)}(n+a_j) \sg_j^{\nu_j}(n+a_j) e(u_j X_{j,y}(n+a_j))\right) d\mbf{u} \\
&\ll x \sum_{1 \leq j \leq k}  \prod_{j' \neq j} \|\hat{w}_{j'}\|_1 \cdot \int_{|u_j| > Z} |\hat{w}_j(u_j)| du_j \ll_{k,\mbf{I},A} x Z^{-A}.
\end{align*}
For each $\nu \in \{0,1\}$, $u \in \mb{R}$ and $n \in \mc{N}$ set 
$$
g_j^{\nu,u}(n) := 1_{\mc{N}}(n) \sg_j^{\nu}(n) e\left(u\frac{\log|\lambda_j(n)|}{\sqrt{c\log\log x}}\right).
$$ 
Clearly, $g_j^{\nu,u}(n)$ is a bounded multiplicative function. Write also
$$
\tilde{W}(\mbf{u}) := \prod_{1 \leq j \leq k} \hat{w}_j(u_j) e\left(\frac{u_j}{\sqrt{4c\log\log x}}\right).
$$
Separating the term $\mbf{\nu} = \mbf{0}$ from the rest, we have
\begin{align}
\tilde{N}_{\mbf{a},y}(x;\mbf{\eps},\tilde{\mbf{I}}) &= 2^{-k}\int_{\|\mbf{u}\|_{\infty} \leq Z} W(\mbf{u}) \sum_{n \leq x} \prod_{1 \leq j \leq k} 1_{\mc{N}^{(\xi)}(2x)}(n+a_j) e(u_jX_{j,y}(n+a_j)) d\mbf{u} \label{eq:nuZero}\\
&+ 2^{-k}\sum_{\mbf{\nu}\in \{0,1\}^k \atop \mbf{\nu} \neq \mbf{0}} \left(\prod_{1 \leq i \leq k} \e_i^{\nu_i}\right) \int_{\|\mbf{u}\|_{\infty} \leq Z} \tilde{W}(\mbf{u}) \left(\sum_{n \leq x} g_1^{\nu_1,u_1}(n+a_1) \cdots g_k^{\nu_k,u_k}(n+a_k)\right) d\mbf{u} \label{eq:nuNonZero}\\
&+ O_{k,\mbf{I},A}(xZ^{-A}) \nonumber =: \Sigma_1 + \Sigma_2 + O_{k,\mbf{I},A}(xZ^{-A}),  \nonumber
\end{align}
where we applied Lemma \ref{lem:approx} (and $Z \leq \log\log X$) to pass from $\mc{N}^{(\xi)}(2x)$ to $\mc{N}$ in the sum over $\mbf{\nu} \neq \mbf{0}$. In the next two subsections, we treat the terms separately, according to whether or not $\mbf{\nu} = \mbf{0}$.
\subsection{The terms with $\mbf{\nu} \neq \mbf{0}$}
Our goal here is to prove the following.
\begin{prop}\label{prop:EllApp}
Let $Y$ be large. Suppose $\mbf{\nu} \neq \mbf{0}$, and let $\|u\|_{\infty} \leq (\log\log X)^{1/3}$. \\
a) Assume that Conjecture \ref{conj:EllIntro} holds. Then
$$
\sum_{n \leq Y} \prod_{1 \leq j \leq k} g_j^{\nu_j,u_j}(n+a_j) = o_{\mbf{a}}(Y).
$$
b) If $1 \leq k \leq 2$ then, unconditionally,
$$
\sum_{n \leq Y} \frac{\prod_{1 \leq j \leq k} g_j^{\nu_j,u_j}(n+a_j)}{n} = o_{\mbf{a}}(\log Y).
$$
c) If $k = 3$ and at least two of the cusp forms $f_1,f_2, f_3$ are the same then, again unconditionally, 
$$
\sum_{n \leq Y} \frac{ g_1^{\nu_1,u_1}(n+a_1)g_2^{\nu_2,u_2}(n+a_2)g_3^{\nu_3,u_3}(n+a_3)}{n} = o_{\mbf{a}}(\log Y).
$$
\end{prop}
Given multiplicative functions $f,g: \mb{N} \ra \mb{U}$ and $Y \geq 2$, we define the \emph{pretentious distance}, in the sense of Granville and Soundararajan, as
$$
\mb{D}(f,g;Y) := \left(\sum_{p \leq Y} \frac{1-\text{Re}(f(p)\bar{g}(p))}{p}\right)^{1/2}.
$$
Variants of Conjecture~\ref{conj:EllIntro} concerning logarithmically-averaged correlations, are known to hold when $k = 2,3$. We may state these as follows.
\begin{thm}[Tao \cite{tao}; Tao-Ter\"{a}v\"{a}inen \cite{TaTe}] \label{thm:logEll} \ \\
a) Let $k = 2$ and $h \geq 1$. Let $g_1,g_2 : \mb{N} \ra \mb{U}$ be multiplicative functions, one of which satisfies \eqref{eq:nonPret}. Then
$$
\sum_{n \leq X} \frac{g_1(n)g_2(n+h)}{n} = o_h(\log X).
$$
b) Let $k = 3$, and let $h_1,h_2$ be distinct positive integers. Let $g_1,g_2,g_3 : \mb{N} \ra \mb{U}$ be multiplicative functions such that for any Dirichlet character $\chi$ one has 
$$
\mb{D}(g_1g_2g_3,\chi;Y)^2 \gg \log\log Y \text{ as $Y \ra \infty$.}
$$ 
Then in this case
$$
\sum_{n \leq x} \frac{g_1(n)g_2(n+h_1)g_3(n+h_2)}{n} = o_{h_1,h_2}(\log x).
$$
\end{thm} 
For the purposes of proving Proposition \ref{prop:EllApp} we establish the following  analogue of Theorem \ref{thm:effST}, establishing a quantitative Sato-Tate type estimate for primes in a progression to a relatively small modulus.
\begin{prop} \label{prop:ThoAPs}
Let $Y$ be large, and let $D \geq 1$. Let $2 \leq q \leq (\log \log Y)^D$ and let $0 \leq a \leq q-1$ be coprime to $q$. Let $J \subset [-2,2]$ be a non-empty interval. Let $f$ be a primitive non-CM cusp form of even weight $r$, squarefree level and trivial nebentypus. Then
$$
\sum_{p \leq Y \atop p \equiv a \pmod{q}} 1_{J}(\lambda_f(p)) = \frac{Li(Y)}{\pi\phi(q)} \left(\int_J \sqrt{1-(t/2)^2} dt + O_{A}\left(Y(\log \log Y)^{-1/2+o(1)}\right)\right),
$$
where $Li(t) := \int_2^t \frac{dv}{\log v}$. 
\end{prop}
\begin{proof}
Write $I := \{\theta \in [0,\pi] : 2\cos \theta \in J\}$. We follow the argument in \cite{Tho}. Using Beurling-Selberg majorants/minorants for $1_I$ (as in \cite{RoTh}), we may select a parameter $M \geq 1$ and functions $F_{I,M}^{\pm}$ such that $F_{I,M}^-(\theta) \leq 1_I(\theta) \leq F_{I,M}^+(\theta)$, and upon writing $F_{I,M}^{\pm}$ in the basis of Chebyshev polynomials of the second kind,
$$
F_{I,M}^{\pm}(\theta) = \sum_{m \in \mb{Z}} \hat{F}_{I,M}^{\pm}(m)U_m(\cos\theta),
$$
we have \\
i) $|\hat{F}_{I,M}^{\pm}(0)-\mu_{ST}(I)| \ll \frac{1}{M+1}$, \\
ii) for $|m| > M$, $\hat{F}_{I,M}^{\pm}(m) \neq 0$, and \\
iii) for $0 < |m| \leq M$ we have 
$$
|\hat{F}_{I,M}^{\pm}(m)| \ll \min\{|I|,1/m\} + \frac{1}{M+1}.
$$ 
Employing these majorants and minorants and their Chebyshev expansions, we obtain
\begin{align*}
&\left|\sum_{p \leq Y \atop p \equiv a \pmod{q}} 1_J(\lambda(p)) - \mu_{ST}(I) \sum_{p \leq Y \atop p \equiv a \pmod{q}} 1\right| \ll \frac{\pi(Y)}{\phi(q)M} + \sum_{1 \leq m \leq M} \frac{1}{m}\left|\sum_{p \leq Y \atop p \equiv a \pmod{q}} U_m(\cos\theta_p) \right| \\
&\leq \frac{\pi(Y)}{\phi(q)M} + \frac{1}{\phi(q)} \sum_{\chi \pmod{q}} \sum_{1 \leq m \leq M} \frac{1}{m} \left|\sum_{p \leq X} U_m(\cos \theta_p) \chi(p)\right|.
\end{align*}
For each $1 \leq m \leq M$ and each non-principal $\chi \pmod{q}$ let $\pi := \text{Sym}^m(f)$, and let $\pi'$ correspond to the $GL_1(\mb{A})$ unitary representation induced by $\chi$ (writing $\mb{A}$ to denote the ad\`{e}les over $\mb{Q}$). By the work of Newton and Thorne \cite{NeTh}, $\pi$ is a cuspidal automorphic representation of $GL_{m+1}(\mb{A})$, and therefore by the Rankin-Selberg theory the convolution $\pi \times \pi'$ is also cuspidal automorphic for $GL_{m+1}(\mb{A})$. Writing the standard automorphic $L$-function of $\pi \times \pi'$ as
$$
L(s,\pi \times \pi') = \sum_{n \geq 1} \frac{a_{\pi \times \pi'}(n)}{n^s} \text{ for $\text{Re}(s) > (r-1)/2+1$,} 
$$
we note that $a_{\pi \times \pi'}(p) = U_m(\cos \theta_p) \chi(p)$ for all primes $p$.  Classical log-free zero density estimates and zero-free regions show that Corollary\footnote{Obviously, if $\chi$ is not quadratic then $\pi'$ is not self-dual. However, since $\pi$ is self-dual we may instead appeal to an analogous, somewhat weaker zero-free region of $L(s,\pi \times \pi')$ due to Brumley \cite[Appendix A]{HuB}. Applying this bound, rather than Corollary 4.2, in the proof of Lemma 6.4 where it is used, we obtain the same result as \cite{Tho}.} 4.2(2) and Proposition 4.4(2) of \cite{Tho} hold for $\pi'$, and Lemma 6.4 there implies that they hold for $\pi$ as well. Then Proposition 5.1(2) of \cite{Tho} shows that
$$
|\sum_{p \leq Y \atop p \nmid q} a_{\pi \times \pi'}(p) \log p| \ll \frac{Y^{\beta_{\chi}}}{\beta_{\chi}} + m^2Y\left(Y^{-c/M^2} + \exp\left(-c\frac{\log Y}{M\log (qC(\pi))}\right) + \exp\left(-c\frac{\sqrt{\log Y}}{M}\right)\right),
$$
where $\beta_{\chi}$ is the possible (necessarily real) Siegel zero of $L(s,\pi \times \pi')$ with $\pi' = \chi$, and $C(\pi)$ is the conductor of $\pi = \text{Sym}^m f$. Equation (6.5) in \cite{Tho} shows that $\log(qC(\pi)) \asymp \log q  + m \log(rm)$. By known bounds for Siegel zeros of automorphic $L$-functions, we have
$$
\beta_{\chi} \leq 1- \frac{c_f}{(qM)^{4M^2}} = 1- c_f \exp\left(-4M^2\log(qM)\right),
$$
for some constant $c_f > 0$ depending only on $f$. Setting $M := C' (\log \log Y)^{1/2}(\log\log \log Y)^{-1}$, where $C'> 0$ is a sufficiently large constant and applying partial summation we obtain, for any $B > 0$, a bound of the shape
$$
\left|\sum_{p \leq Y \atop p \nmid q} U_m(\cos \theta_p) \chi(p)\right| \ll_{A,B} m^2Y (\log Y)^{-B}.
$$
It follows that if $q \leq (\log \log Y)^D$ we can choose $B = B(D)$ such that
$$
\sum_{p \leq Y \atop p \equiv a \pmod{q}} 1_J(\lambda(p)) = \frac{1}{\phi(q)}\left(\mu_{ST}(I)\text{Li}(Y) + O_D\left(\frac{Y}{(\log Y)(\log\log Y)^{1/2-o(1)}} \right)\right),
$$
and the claim follows.
\end{proof}
In order to apply Elliott's conjecture we need to establish that at least one of the functions $g_j^{\nu_j,u_j}$ satisfies \eqref{eq:nonPret}. In fact, we show the stronger estimate. 

\begin{lem} \label{lem:nonPret}
Let $Y$ be large. For any $\mbf{\nu} \neq \mbf{0}$ and $\mbf{u} \in \mb{R}^k$ with $\|\mbf{u}\|_{\infty} \leq (\log\log Y)^{1/3}$ there is a $1 \leq j_0 \leq k$ such that for any $A > 1$, if $Q := (\log\log Y)^A$ then
$$
\min_{|t| \leq Y} \min_{\chi \pmod{q} \atop q \leq Q} \mb{D}(g_{j_0}^{\nu_{j_0},u_{j_0}},\chi n^{it}; Y)^2 \gg \log\log Y.
$$
\end{lem}
\begin{proof}
Since $\mbf{\nu} \neq \mbf{0}$ we may fix an index $1 \leq j_0 \leq k$ for which $\nu_{j_0} = 1$, and as such 
$$
g_{j_0}^{\nu_{j_0},u_{j_0}}(n) = \sg_{j_0}(n)1_{\mc{N}}(n) e\left(u_{j_0}\frac{\log|\lambda_{j_0}(n)|}{c\sqrt{\log\log X}}\right) \text{ for all $n \in \mb{N}$.}
$$
For convenience, write $g_{j_0} = g_{j_0}^{\nu_{j_0},u_{j_0}}$ for the remainder of the proof. Also, let $(t_0,\psi)$ be the pair consisting of $|t_0| \leq Y$ and $\psi$ a Dirichlet character to a modulus $\leq Q$ that minimizes the expression on the LHS; we may assume that $\psi$ is primitive and we denote its conductor by $q$. We note first of all that whenever $|\lambda(p)| \geq (\log\log Y)^{-1}$ we have 
$$
g_{j_0}(p) = \sg_{j_0}(p) + O(|u_{j_0}|(\log\log \log Y)(\log\log Y)^{-1/2}),
$$ 
so in view of the bound on $\|\mbf{u}\|_{\infty}$ and Lemma \ref{lem:LRSSharp} we obtain
\begin{align} \label{eq:remU}
&\min_{|t| \leq Y} \min_{\chi \pmod{q} \atop q \leq Q} \mb{D}(g_{j_0},\chi n^{it};Y)^2 \nonumber\\
&\geq \sum_{p \leq Y \atop |\lambda_{j_0}(p)| \geq (\log\log Y)^{-1}} \frac{1-\text{Re}(\sg_{j_0}(p) \bar{\psi}(p)p^{-it_0})}{p} + O(\|\mbf{u}\|_{\infty} (\log\log \log Y) (\log\log Y)^{-1/2})  \nonumber \\
&= \sum_{p \leq Y} \frac{1-\text{Re}(\sg_{j_0}(p)\bar{\psi}(p)p^{-it_0})}{p} + o(\log\log Y).
\end{align}
Appendix C of \cite{MRT} shows that $\mb{D}(\sg_{j_0},\psi n^{it_0};Y)^2 \geq \frac{1}{16} \log\log Y$ provided that either $\psi^2$ is non-principal or $1 \leq |t_0| \leq Y$, and the claim then follows in either of these cases. \\
Thus, we may assume that $\psi$ is quadratic and $|t_0| \leq 1$.
It suffices to prove now that
$$
\sum_{p \leq Y} \frac{\sg_{j_0}(p)\psi(p)}{p^{1+it_0}} = o(\log\log Y),
$$
as $Y \ra \infty$. Put $Z := \exp((\log Y)^{1/\log\log\log Y})$, so that if $Y$ is large enough and $A > 1$ then $q \leq (\log\log Z)^{A+1}$. By partial summation and $|t_0| \leq 1$,
\begin{align*}
\left|\sum_{p \leq Y} \frac{\sg_{j_0}(p)\psi(p)}{p^{1+it_0}}\right| &= |1+it_0| \left|\int_Z^Y \left(\sum_{p \leq u} \sg_{j_0}(p) \psi(p)\right)\frac{du}{u^{2+it_0}}\right| + O(\log\log Z) \\
&\ll \int_Z^Y \left|\sum_{p \leq u} \sg_{j_0}(p)\psi(p)\right| \frac{du}{u^2} + o(\log\log Y).
\end{align*}
It is therefore enough to show that if $u \geq Z$ then 
$$
\left|\sum_{p \leq u} \sg_{j_0}(p)\psi(p)\right| \ll \frac{u}{(\log u)(\log\log u)^{1/2-o(1)}}.
$$
The sum in brackets can be written as
\begin{equation} \label{eq:STExpand}
\asum_{a \pmod{q}} \psi(a) \left(\sum_{p \leq u \atop p \equiv a \pmod{q}} 1_{\lambda_{j_0}(p) > 0} - \sum_{p \leq u \atop p \equiv a \pmod{q}} 1_{\lambda_{j_0}(p) < 0}\right).
\end{equation}
Applying Proposition \ref{prop:ThoAPs} with $J = [-2,0]$ and $J = [0,2]$ (and subtracting the contribution from $J = \{0\}$ in each), we have
$$
\sum_{p \leq u \atop p \equiv a \pmod{q}} 1_{J}(\lambda_{j_0}(p)) = \left(\frac{1}{\pi}+O\left((\log \log u)^{-1/2+o(1)}\right)\right) \left(\int_J \sqrt{1-(t/2)^2}dt \right) \cdot \frac{u}{\phi(q)\log u}.
$$
Since the Sato-Tate distribution is symmetric, the bracketed expression in \eqref{eq:STExpand} is
\begin{align*}
&\frac{1}{\pi\phi(q)} \frac{u}{\log u} \left(\int_{-2}^0 \sqrt{1-(t/2)^2}dt - \int_0^2 \sqrt{1-(t/2)^2} dt  + O(u(\log \log u)^{-1/2+o(1)})\right) \\
&= O\left(\frac{u}{\phi(q)(\log u)(\log\log u)^{1/2-o(1)}}\right),
\end{align*}
which proves the desired estimate upon summing over all $a \pmod{q}$ coprime to $q$.
\end{proof}
\begin{proof}[Proof of Proposition \ref{prop:EllApp}]
a) This part follows immediately, assuming Conjecture \ref{conj:EllIntro}, from the conclusion of Lemma \ref{lem:nonPret}. \\
Note that (up to relabeling) it suffices to show the case where $\nu_j = 1$ for $1 \leq j \leq k$, for all $1 \leq k \leq 3$, in parts b) and c). \\
When $k = 1$, Theorem 1.6 of \cite{LamMan} (with $T = 1$) shows that
$$
\sum_{n \leq x} \frac{g_1^{1,u_1}(n+a_1)}{n} = \sum_{n \leq x} \frac{g_1^{1,u_1}(n)}{n} + O_{\mbf{a}}(1) \ll_{\mbf{a}} 1+ (\log x)\exp\left(-\min_{|t| \leq 1} \mb{D}(g_1,n^{it};x)^2\right) = o(\log x),
$$
by Lemma \ref{lem:nonPret} (here we may take $Q = 1$). \\
When $k = 2$, the result follows directly upon combining Theorem \ref{thm:logEll} a) with Lemma \ref{lem:nonPret}. \\
When $k = 3$ we may observe (as we did in the proof of Lemma \ref{lem:nonPret}) that on primes $p$ for which $|\lambda_j(p)| \geq (\log\log X)^{-1}$ for all $j = 1,2,3$ we have
$$
g_1^{1,u_1}g_2^{1,u_2}g_3^{1,u_3}(p) = \sg_1(p)\sg_2(p)\sg_3(p)+ O\left(\|\mbf{u}\|_{\infty} \frac{\log\log\log X}{\sqrt{\log\log X}}\right).
$$
In particular, if at least two of the forms, say $f_1$ and $f_2$, are the same, then $\sg_1(p) = \sg_2(p)$ for all $p$, and thus $\sg_1(p)\sg_2(p)\sg_3(p) = \sg_3(p)$ for all primes $p$ such that $|\lambda(p)| \geq (\log\log x)^{-1}$, say. We obtain, in analogy to \eqref{eq:remU},
$$
\mb{D}(g_1^{1,u_1}g_2^{1,u_2}g_3^{1,u_3}, \chi; Y)^2 \geq \sum_{p \leq Y} \frac{1-\text{Re}(\sg_3(p)\bar{\chi}(p))}{p} + o(\log\log Y),
$$
for any Dirichlet character $\chi$, as $Y \ra \infty$. It follows from the remainder of the proof of Lemma \ref{lem:nonPret} that $\mb{D}(\sg_3,\chi;X)^2 \gg \log\log X$, and so by Theorem \ref{thm:logEll} b) we may conclude that
$$
\sum_{n \leq X} \frac{g_1^{1,u_1}(n+a_1)g_2^{1,u_2}(n+a_2)g_3^{1,u_3}(n+a_3)}{n} = o(\log X),
$$
and the proof is complete.
\end{proof}

\subsection{The term $\mbf{\nu} = \mbf{0}$}
To treat the contribution from $\mbf{\nu} = \mbf{0}$ we seek to estimate the quantity
$$
\phi_{\mbf{a},y}(\mbf{u}) := \sum_{n \leq x} \prod_{1 \leq j \leq k} 1_{\mc{N}^{(\xi)}(2x)}(n+a_j) e\left(u_jX_{j,y}(n+a_j)\right).
$$
A moment's reflection allows us to identify this as the unnormalized characteristic function (in the sense of probability theory) of the distribution function
$$
\mbf{t} \mapsto X^{-1}|\{n \leq x : n+a_j \in \mc{N}^{(\xi)}(2x), X_{j,y}(n+a_j) \leq t_j \text{ for all } 1 \leq j \leq k\}|.
$$
We establish the following asymptotic formula for $\phi_{\mbf{a},y}$.
\begin{prop}\label{prop:charfn}
Suppose $\|\mbf{u}\|_{\infty} \leq Z$. Then
$$
\phi_{\mbf{a},y}(\mbf{u}) = \exp\left(-2\pi^2\sum_{1 \leq j \leq k} u_j^2\right) + O_k\left(Z^{k-1}2^{-Z^2} + \frac{(\log\log \log x)^{O(1)}}{\log\log x} e^{4\pi^2 k Z^2}\right).
$$
\end{prop}
\begin{proof}
Let $N \geq 1$ be a parameter to be chosen. Using the identity
$$
e(t) = \sum_{0 \leq j \leq 2N-1} \frac{(2\pi i t)^j}{j!} + O\left(\frac{(2\pi t)^{2N}}{(2N!)}\right),
$$
valid for any $t \in \mb{R}$, we obtain
\begin{align*}
\phi_{\mbf{a},y}(\mbf{u}) &= \sum_{\mbf{l} \in (\mb{N} \cup \{0\})^k \atop \|\mbf{l}\|_{\infty} \leq 2N-1} \left(\prod_{1 \leq j \leq k} \frac{(2\pi i u_j)^{l_j}}{l_j!}\right) \sum_{n \leq x} \prod_{1 \leq j \leq k} 1_{\mc{N}^{(\xi)}(2x)}(n+a_j) X_{j,y}(n+a_j)^{l_j} \\
&+O\left(\sum_{\mbf{l} \in (\mb{N} \cup \{0\})^k \atop \|\mbf{l}\|_{\infty} = 2N}\frac{(2\pi |u_j|)^{l_j}}{l_j!} \sum_{n\leq x} \prod_{1 \leq j \leq k} 1_{\mc{N}^{(\xi)}(2x)}(n+a_j) |X_{j,y}(n+a_j)|^{l_j}\right) =: T_1 + T_2.
\end{align*}
We recognize the inner sum in the main term $T_1$ as $\tilde{M}_{\mbf{f},\mbf{l}}(x) = M_{\mbf{f},\mbf{l}}(x)\prod_{1 \leq j \leq k} \sg_j(y)^{-l_j}$, and by Proposition \ref{prop:Moments} we get
\begin{align*}
T_1 &= X\sum_{\mbf{l} \in (\mb{N} \cup \{0\})^k \atop \|\mbf{l}\|_{\infty} \leq 2N-1} \prod_{1 \leq j \leq k} \frac{(2\pi i u_j)^{l_j}}{l_j!} \cdot \frac{\Gamma(l_j+1)}{2^{l_j/2}\Gamma(l_j/2+1)}\left(1_{2|l_j}  + O_k\left(\frac{l_j^3(\log\log\log x)^{O(1)}}{\log\log x}\right)\right) \\
&= X\left(\sum_{\mbf{r} \in (\mb{N}\cup \{0\})^k \atop \|\mbf{r}\|_{\infty} \leq N-1} \frac{(2\pi i u_j)^{2r_j}}{2^{r_j}r_j!} + O_k\left(\frac{(\log\log\log x)^{O(1)}}{\log\log x} \sum_{\mbf{l} \in (\mb{N}\cup \{0\})^k \atop \|\mbf{l}\|_{\infty} \leq 2N-1} \frac{(2\pi |u_j|)^{l_j}l_j^3}{2^{l_j/2}\Gamma(l_j/2+1)}\right)\right) \\
&= X\left(\exp\left(-2\pi^2 \sum_{1 \leq j \leq k} u_j^2\right) + O_k\left(\frac{Z^{2N}}{N!}+\frac{N^3(\log\log\log x)^{O(1)}}{\log\log x} e^{4\pi^2 k Z^2}\right)\right).
\end{align*} 
Next, we treat $T_2$. For $\mbf{l} \in (\mb{N} \cup \{0\})^k$, define associated vectors $\mbf{l}^{\pm} = (l_1^{\pm},\ldots,l_k^{\pm})$ via
\begin{align*}
l_j^+ := \begin{cases} l_j &\text{ if $2|l_j$,} \\ l_j+1 &\text{ if $2 \nmid l_j$,} \end{cases} \ \ \ \ \ \ \ \ 
l_j^- := \begin{cases} l_j &\text{ if $2|l_j$,} \\ l_j-1 &\text{ if $2 \nmid l_j$,} \end{cases}
\end{align*}
for all $1 \leq j \leq k$. Since $2\mbf{l} = \mbf{l}^+ + \mbf{l}^-$, we obtain, by Cauchy-Schwarz
\begin{align*}
T_2 &\ll \sum_{\mbf{l} \in (\mb{N} \cup \{0\})^k \atop \|\mbf{l}\|_{\infty} = 2N}\left(\prod_{1 \leq j \leq k} \frac{(2\pi |u_j|)^{l_j}}{l_j!}\right) \sum_{n \leq x}  \left(\prod_{1 \leq j \leq k} \prod_{v_j \in \{l_j^+,l_j^-\}} 1_{\mc{N}^{(\xi)}(2x)}(n+a_j)|X_{j,y}(n+a_j)|^{v_j/2}\right) \\
&\ll \sum_{\mbf{l} \in (\mb{N} \cup \{0\})^k \atop \|\mbf{l}\|_{\infty} = 2N}\left(\prod_{1 \leq j \leq k} \frac{(2\pi |u_j|)^{l_j}}{l_j!}\right)\tilde{M}_{\mbf{f},\mbf{l}^+}(x)^{1/2}\tilde{M}_{\mbf{f},\mbf{l}^-}(x)^{1/2}.
\end{align*}
Since $2|l_j^{\pm}$ for all $1 \leq j \leq k$, Proposition \ref{prop:Moments} once again implies that
\begin{align*}
T_2 &\ll X\sum_{\mbf{l} \in (\mb{N} \cup\{0\})^k \atop \|\mbf{l}\|_{\infty} = 2N} \prod_{1 \leq j \leq k} \frac{(2\pi |u_j|)^{l_j}}{l_j!}\cdot \left(\frac{l_j^+!l_j^-!}{2^{l_j}(l_j^+/2)!(l_j^-/2)!}\right)^{1/2} \\
&\ll X\sum_{\mbf{l} \in (\mb{N} \cup\{0\})^k \atop \|\mbf{l}\|_{\infty} = 2N} \prod_{1 \leq j \leq k} \frac{(2\pi |u_j|)^{l_j}}{2^{l_j/2}\Gamma(l_j/2+1)},
\end{align*}
using the fact that the Gamma function is log convex in the last line. Bounding this last expression trivially, we obtain
\begin{align*}
T_2 &\ll X\sum_{1 \leq j \leq k} \frac{(2\pi |u_j|)^{2N}}{2^{N}N!} \sum_{\mbf{r} \in (\mb{N} \cup \{0\})^{k-1} \atop \|\mbf{r}\|_{\infty} \leq N}\prod_{1 \leq j' \leq k \atop j' \neq j}  \frac{(2\pi |u_{j'}|)^{2r_{j'}}}{2^{r_{j'}}}\left(\frac{1}{\Gamma(r_{j'}+1)} + \frac{|u_{j'}|}{2\Gamma(r_{j'}+3/2)}\right) \\
&\ll X\frac{(2\pi)^{2N}Z^{2N+k}}{2^N N!} \sum_{1 \leq j \leq k} \prod_{1 \leq j' \leq k \atop j' \neq j} \sum_{0 \leq r_{j'}\leq N} \frac{(2\pi^2|u_{j'}|^2)^{r_{j'}}}{r_{j'}!} \\
&\ll_k X\frac{(2\pi^2)^{N}Z^{2N+k}}{N!} e^{2\pi^2 (k-1)Z^2}.
\end{align*}
We conclude that
$$
\phi_{\mbf{a},y}(\mbf{u}) = X\left(\exp\left(-2\pi^2\sum_{1 \leq j \leq k} u_j^2\right) + O_k\left(\frac{(2\pi^2)^N Z^{2N+k}}{N!} e^{2\pi^2(k-1)Z^2} + \frac{N^3(\log\log\log x)^{O(1)}}{\log\log x} e^{4 \pi^2 kZ^2}\right)\right).
$$
We select $N = \left \lceil 4\pi^2 e^{k} Z^2\right \rceil \ll_k (\log\log\log x)^{O(1)}$. By Stirling's formula,
\begin{align*}
\phi_{\mbf{a},y}(\mbf{u}) &= X\left(\exp\left(-2\pi^2\sum_{1 \leq j \leq k} u_j^2\right) + O_k\left(\frac{Z^k}{\sqrt{N}} \left(\frac{2\pi^2 Z^2 e^k}{N}\right)^N + \frac{(\log\log \log x)^{O_k(1)}}{\log\log x} e^{4\pi^2 k Z^2}\right)\right) \\
&= X\left(\exp\left(-2\pi^2\sum_{1 \leq j \leq k} u_j^2\right) + O_k\left(Z^{k-1}2^{-Z^2} + \frac{(\log\log \log x)^{O_k(1)}}{\log\log x} e^{4\pi^2 k Z^2}\right)\right),
\end{align*}
and the claim follows.
\end{proof}
\begin{proof}[Proof of Proposition \ref{prop:EKwithSigns}a)]
Assume Conjecture \ref{conj:EllIntro}. 
We choose $Z := 2(\log\log\log x)^{1/3}$. In light of \eqref{eq:nuZero} and \eqref{eq:nuNonZero}, it suffices to show the estimates
\begin{align*}
\Sigma_1 & = 2^{-k}\int_{\|\mbf{u}\|_{\infty} \leq Z} W(\mbf{u}) \sum_{n \leq x} \prod_{1 \leq j \leq k} 1_{\mc{N}^{(\xi)}(2x)}(n+a_j) e(u_jX_{j,y}(n+a_j)) d\mbf{u} \\
&= X\left( 2^{-k} \cdot (2\pi)^{-k/2}\int_{\mc{B}(\mbf{I})} e^{-\frac{1}{2}(x_1^2 + \cdots + x_k^2)} d\mbf{x} + o(1)\right) + O_{k,\mbf{a}}(\delta x)
\end{align*}
as well as 
\begin{align*}
\Sigma_2 &= 2^{-k}\sum_{\mbf{\nu}\in \{0,1\}^k \atop \mbf{\nu} \neq \mbf{0}} \left(\prod_{1 \leq i \leq k} \e_i^{\nu_i}\right) \int_{\|\mbf{u}\|_{\infty} \leq Z} \tilde{W}(\mbf{u}) \left(\sum_{n \leq x} g_1^{\nu_1,u_1}(n+a_1) \cdots g_k^{\nu_k,u_k}(n+a_k)\right) d\mbf{u} \\
&= o_{k,\mbf{a},A,\mbf{I}}(x),
\end{align*}
after which point we will have
$$
N_{\mbf{a},y}(x;\mbf{\eps},\tilde{\mbf{I}}) = \tilde{N}_{\mbf{a},y}(x;\mbf{\eps},\tilde{\mbf{I}}) - O_k(\delta x) =  X\left(2^{-k}\cdot (2\pi)^{-k/2} \int_{\mc{B}(\mbf{I})} e^{-\frac{1}{2}(x_1^2+\cdots + x_k^2)} d\mbf{x} - O_{k,\mbf{a}}(\delta)\right),
$$
given that $\mc{B}(\mbf{I}) \bk \mc{B}(\tilde{\mbf{I}})$ has measure $O_k(\delta)$. Using \eqref{eq:smoothSign} and the line that follows it, we thus get (with $\tilde{\mbf{I}} = \mbf{I}^{(1)}, \mbf{I}^{(2)}$)
\begin{align*}
&X\left(2^{-k}\cdot (2\pi)^{-k/2} \int_{\mc{B}(\mbf{I})} e^{-\frac{1}{2}(x_1^2+\cdots + x_k^2)} d\mbf{x} - O_{k,\mbf{a}}(\delta)\right) \leq N_{\mbf{a}}(x; \mbf{\eps},\mbf{I}) + o(x) \\
&\leq X\left(2^{-k}\cdot (2\pi)^{-k/2} \int_{\mc{B}(\mbf{I})} e^{-\frac{1}{2}(x_1^2+\cdots + x_k^2)} d\mbf{x} + O_{k,\mbf{a}}(\delta)\right),
\end{align*}
whenever $x$ is sufficiently large in terms of $\delta$. Since $\delta$ is arbitrary small, the claim will follow. \\
We begin with $\Sigma_2$. By the triangle inequality, we of course have
\begin{align*}
|\Sigma_2| \leq 2^{-k}\|\tilde{W}\|_1 \sum_{\mbf{\nu} \in \{0,1\}^k \atop \mbf{\nu} \neq \mbf{0}} \max_{\|\mbf{u}\|_{\infty} \leq Z} \left|\sum_{n \leq x} g_1^{\nu_1,u_1}(n+a_1) \cdots g_k^{\nu_k,u_k}(n+a_k)\right|.
\end{align*}
By \eqref{eq:FTbd}, we have $\|\tilde{W}\|_1 = \prod_{1 \leq j \leq k} \|\hat{w}_j\|_1 \ll_{k,\mbf{I}} 1$. Now, since the map 
$$
\mbf{u} \mapsto \sum_{n \leq X} \prod_{1 \leq j \leq k} g_j^{\nu_j,u_j}(n+a_j)
$$
is continuous, its maximum over the compact region $[-Z,Z]^k$ is attained at some point $\mbf{u}_0(\mbf{\nu}) = (u_{0,1}(\mbf{\nu}),\ldots,u_{0,k}(\mbf{\nu}))$, say, for each $\mbf{\nu} \neq \mbf{0}$. By Proposition \ref{prop:EllApp} a), we thus have
\begin{align*}
2^{-k}\sum_{\mbf{\nu} \in \{0,1\}^k \atop \mbf{\nu} \neq \mbf{0}} \max_{\|\mbf{u}\|_{\infty} \leq Z} \left|\sum_{n \leq x} \prod_{1 \leq j \leq k} g_j^{\nu_j,u_j}(n+a_j)\right| = 2^{-k}\sum_{\mbf{\nu} \in \{0,1\}^k \atop \mbf{\nu} \neq \mbf{0}} \left|\sum_{n \leq x} \prod_{1 \leq j \leq k} g_j^{\nu_j,u_{0,j}(\mbf{\nu})}(n+a_j)\right| = o_{\mbf{a}}(X).
\end{align*}
It follows therefore that
$$
\Sigma_2 = o_{k,\mbf{a},\mbf{I}}(x),
$$
as claimed. \\
Next, we estimate $\Sigma_1$. In light of our choice of $Z$, Lemma \ref{prop:charfn} yields
\begin{align*}
\Sigma_1 &= 2^{-k}\int_{\|\mbf{u}\|_{\infty} \leq Z} W(\mbf{u}) \phi_{\mbf{a},y}(\mbf{u}) d\mbf{u} = \left(2^{-k} \int_{\|\mbf{u}\|_{\infty} \leq Z} \prod_{1 \leq j \leq k} \hat{w}_j(u_j) e^{-2\pi^2 u_j^2} d\mbf{u}\right) X + o_{k,A,\mbf{a}}(x) \\
&= \left(2^{-k} \prod_{1 \leq j \leq k} \left(\int_{\mb{R}} \hat{w}_j(u_j) e^{-2\pi^2u_j^2} du_j\right)\right) X + o_{k,A,\mbf{a}}(x).
\end{align*}
By Plancherel's theorem, we obtain
$$
\Sigma_1 = \left(2^{-k} (2\pi)^{-k/2} \prod_{1 \leq j \leq k} \int_{\mb{R}} w_j(y_j) e^{-y_j^2/2} dy_j\right) X + o_{k,A,\mbf{a}}(x).
$$
Since for each $1 \leq j \leq k$ we have that $w_j \geq 0$, $w_j \equiv 1$ on $I_j$ and it is supported on $I_{j,\delta}$, it follows that
\begin{align*}
\Sigma_1 &= \left(2^{-k} (2\pi)^{-k/2} \prod_{1 \leq j \leq k} \int_{I_j} e^{-y_j^2/2} dy_j\right)X + O_{k,\mbf{a}}\left(x\sum_{1 \leq j \leq k} \int_{I_{j,\delta} \bk I_j} e^{-\frac{1}{2}y_j^2} dy_j\right) +  o_{k,A,\mbf{a}}(x) \\
&= \left(2^{-k} (2\pi)^{-k/2} \prod_{1 \leq j \leq k} \int_{I_j} e^{-y_j^2/2} dy_j\right)X + O_{k,\mbf{a}}(\delta x),
\end{align*}
as required.
\end{proof}
\begin{proof}[Proof of Proposition \ref{prop:EKwithSigns} b)]
The proof of part b) is similar to that of a), so we merely highlight the main differences. \\
As with $N_{\mbf{a}}(x;\mbf{\eps},\mbf{I})$, instead of the (logarithmically-averaged) sum over $S_{\mbf{a}}(\mbf{\eps},\tilde{\mbf{I}})$ it suffices to consider
$$
\frac{1}{\log x} \sum_{n \leq x} \frac{1}{n}\prod_{1 \leq j \leq k} 1_{\mc{N}^{(\xi)}(2x)}(n+a_j) 1_{\sg_j(n+a_j) = \e_j} 1_{\tilde{I}_j}(X_{j,y}(n+a_j)).
$$
We smooth out the indicator function of $1_{\tilde{I}_j}$ as above (note that the logarithmically-averaged variants of the asymptotic formulae from the previous section also hold by partial summation), so once again it suffices to consider
$$
\frac{1}{\log x} \sum_{n \leq x} \frac{1}{n}\prod_{1 \leq j \leq k} 1_{\mc{N}^{(\xi)}(2x)}(n+a_j) 1_{\sg_j(n+a_j)} w_j(n+a_j),
$$
with the same choice of smooth functions $w_1,\ldots,w_k$ (and a suitable, small parameter $\delta > 0$) from the previous subsection. \\
Continuing the argument of the previous subsection, using Fourier inversion and \eqref{eq:pmIden}, we conclude that it suffices to check that
\begin{align*}
&\int_{\|\mbf{u}\|_{\infty} \leq Z} W(\mbf{u}) \left(\sum_{n \leq x} \frac{1}{n}\prod_{1 \leq j \leq k} 1_{\mc{N}^{(\xi)}(2x)} (n+a_j) e(u_jX_{j,y}(n+a_j))\right) d\mbf{u} \\
&= (2\pi)^{-k/2} \int_{\mc{B}(\mbf{I})} e^{-\frac{1}{2}(x_1^2+\cdots + x_k^2)} d\mbf{x} + O_{k,\mbf{a}}(\delta x), \\
&\sum_{\mbf{\nu} \in \{0,1\}^k \atop \mbf{\nu} \neq \mbf{0}} \left|\int_{\|\mbf{u}\|_{\infty} \leq Z} \tilde{W}(\mbf{u}) \left(\sum_{n \leq x} \frac{1}{n}\prod_{1 \leq j \leq k} g_j^{\nu_j,u_j}(n+a_j)\right) d\mbf{u}\right| = o_{k,A,\mbf{a},\mbf{I}}(\log x)
\end{align*}
The proof of the first estimate is exactly analogous to the treatment of $\Sigma_1$ in the deduction of part a) (using logarithmically-averaged variants of the moment calculations from the previous section). To treat $\Sigma_2$ we use Proposition \ref{prop:EllApp} b) and c), depending on whether $k = 2$ or $k = 3$, instead of a). The proof then concludes as above.
\end{proof}

\section{Proof of Theorems \ref{thm:genCusp}, \ref{thm:condSigns} and \ref{thm:uncondSigns}} \label{sec:deduc}
\begin{proof}[Proof of Theorem \ref{thm:genCusp}]
Define 
\begin{align*}
T_{\mbf{a}}(x) := |\{n \leq x : n+a_j \in \mc{N}, |b_j(n+a_j)| < |b_{j+1}(n+a_{j+1})| \text{ for all } 1 \leq j \leq k-1\}|.
\end{align*}
Clearly, we have 
$$
|\{n \leq x : |b_1(n+a_1)| < \cdots < |b_k(n+a_k)|\}| \geq T_{\mbf{a}}(x),
$$
and we will show that 
$$
T_{\mbf{a}}(x) = (1/k! + o(1)) X \text{ as $x \ra \infty$,}
$$
so that as $X \gg_{\mbf{a}} x$ (as in Remark \ref{rem:NjPD}), the claim will follow. \\
Note that if $x$ is large, $\sqrt{x} < n \leq x$ and $|\lambda_j(n+a_j)| < |\lambda_{j+1}(n+a_{j+1})|\left(1-\frac{1}{\sqrt{n}}\right)$ for all $1 \leq j \leq k-1$ then $n$ is counted by $T_{\mbf{a}}(x)$: indeed, we then have
\begin{align*}
|b_j(n+a_j)| &< |b_{j+1}(n+a_{j+1})| \left(1+ \frac{a_j-a_{j+1}}{n+a_{j+1}}\right)^{(m-1)/2}\left(1-\frac{1}{\sqrt{n}}\right) \\
&= |b_{j+1}(n+a_{j+1})|\left(1-\frac{1}{\sqrt{n}} + O_{m,\mbf{a}}(1/n)\right) < |b_{j+1}(n+a_{j+1})|,
\end{align*}
for all $1 \leq j \leq k$. \\
For each $n \in \mc{N}$ and $1 \leq j \leq k$ set
$$
\tilde{X}_j(n) := \frac{\log|\lambda_j(n)|+\frac{1}{2}\log\log x}{\sqrt{c\log\log x}}.
$$
Then by the above and monotonicity, we have
\begin{align*}
&T_{\mbf{a}}(x) +O(\sqrt{x}) \\
&\geq |\{\sqrt{x} < n \leq x : n+a_j \in \mc{N}, |\lambda_j(n+a_j)| < |\lambda_{j+1}(n+a_{j+1})|\left(1-\frac{1}{\sqrt{n}}\right) \forall \  1 \leq j \leq k-1\}|\\
&= |\{\sqrt{x} < n \leq x : n+a_j \in \mc{N}, \tilde{X}_j(n+a_j) < \tilde{X}_{j+1}(n+a_{j+1}) + \frac{\log(1-1/\sqrt{n})}{\sqrt{c\log\log x}} \ \forall \ 1 \leq j \leq k-1\}|.
\end{align*}
Moreover, whenever $\sqrt{x} < n \leq x$, $n+a_j \in \mc{N}$ for all $j$ and $x$ is large enough we have
$$
\max_{1 \leq j \leq k} |X_j(n+a_j)-\tilde{X}_j(n+a_j)| \ll \frac{1}{\sqrt{\log\log x}}.
$$
It follows that if $\sqrt{x} < n \leq x$, $n+a_j \in \mc{N}$ and
$$
X_j(n+a_j) < X_{j+1}(n+a_{j+1}) - \frac{1}{(\log\log x)^{1/3}}
$$
for all $1 \leq j \leq k-1$ then for large enough $x$,
\begin{align*}
\tilde{X}_j(n+a_j) &< X_j(n+a_j) + O(1/\sqrt{\log\log x}) < X_{j+1}(n+a_{j+1}) - \frac{1}{2(\log\log x)^{1/3}} \\
&< \tilde{X}_{j+1}(n+a_{j+1}) -\frac{1}{4(\log\log x)^{1/3}} < \tilde{X}_{j+1}(n+a_{j+1}) + \frac{\log(1-1/\sqrt{n})}{\sqrt{\log\log x}},
\end{align*}
for all $1 \leq j \leq k-1$. We therefore find that
\begin{align}
&T_{\mbf{a}}(x) +O(\sqrt{x}) \nonumber\\
&\geq |\{\sqrt{x} < n \leq x : n+a_j \in \mc{N}, X_j(n+a_j) < X_{j+1}(n+a_{j+1}) -\frac{1}{(\log\log x)^{1/3}} \forall \ 1 \leq j \leq k-1\}|. \label{eq:lowBdIneq}
\end{align}
We will bound this quantity from below, but we first make some preliminary remarks.  \\
Put $\mc{B} := \{\mbf{x} \in \mb{R}^k : x_1 < \cdots < x_k\}$.  Note that the set
$$
\{\mbf{x} \in \mb{R}^k : \exists i \neq j \text{ such that } x_i = x_j\}
$$
has Lebesgue measure 0, so by symmetry, we easily find that
$$
(2\pi)^{-k/2} \int_{\mc{B}} e^{-\frac{1}{2}(x_1^2+\cdots + x_k^2)}dx_1\cdots dx_k = \frac{1}{k!}.
$$
In what follows we show that the expression in \eqref{eq:lowBdIneq} is well-approximated by this integral. \\
Given a large integer $Z \geq 1$, define 
$$
\mc{B}(Z) := \{\mbf{x} \in \mb{R}^k : \|\mbf{x}\|_{\infty} \leq Z, x_1 < \cdots < x_k\}.
$$
Clearly, by the union bound we have
$$
\int_{\mc{B} \bk \mc{B}(Z)} e^{-\frac{1}{2}\mbf{x} \cdot \mbf{x}} d\mbf{x} \ll \sum_{1 \leq j \leq k} \left(\prod_{1 \leq l \leq k \atop l \neq j} \int_{\mb{R}} e^{-\frac{1}{2}x_l^2} dx_l\right)\int_Z^{\infty} e^{-\frac{1}{2}x_j^2} dx_j  \ll_k e^{-Z^2/4},
$$
so it suffices to consider the integral over the bounded set $\mc{B}(Z)$.
Given a second integer $N \geq 2$, for each $0 \leq l \leq 2ZN-1$ set 
$$
y_l := -Z + l/N \text{ and } I_l := [y_l,y_{l+1}-(ZN)^{-2k}].
$$ 
Given a $k$-tuple $\mbf{l} = (l_1,\ldots,l_k) \in \{0,\ldots,2ZN-1\}^k$, define boxes
$$
B_{N,Z}(\mbf{l}) := \prod_{1 \leq i \leq k} I_{l_i},
$$
and finally define
$$
\mc{B}_{N}(Z) := \bigcup_{\substack{\mbf{m} \in \{0,\ldots,2ZN-1\}^k \\ m_j < m_{j+1} \\ \forall 1 \leq j \leq k-1}} B_{N,Z}(\mbf{m}),
$$
the union here being disjoint. We easily find that $\mc{B}_N(Z) \subseteq \mc{B}(Z)$, since any $\mbf{x} \in \mc{B}_{N}(Z)$ must lie in some $B_{N,Z}(\mbf{m})$, where $\mbf{m} \in \{0,\ldots,2ZN-1\}$ satisfies $m_j < m_{j+1}$ for all $1 \leq j \leq k-1$, and thus
$$
x_j \leq y_{m_j+1}-(N Z)^{-2k} < y_{m_{j+1}} \leq x_{j+1},
$$ 
for all $1 \leq j \leq k-1$. \\
Conversely, we claim that $\mc{B}(Z) \bk \mc{B}_N(Z)$ has small Lebesgue measure. To see this, suppose $\mbf{x} \in \mc{B}(Z)$. Then for each $1 \leq j \leq k$ there is $m_j \in \{0,\ldots,2ZN-1\}$ such that $x_j \in [y_{m_j},y_{m_j+1}]$ and $m_j \leq m_{j+1}$ for all $1 \leq j \leq k-1$. Thus, if $\mbf{x} \in \mc{B}(Z)\bk \mc{B}_{N}(Z)$ then either: 
\begin{enumerate}[label = (\roman*)]
\item for at least one $1 \leq j_0 \leq k-1$ we have $m_{j_0} = m_{j_0+1}$, and thus $|x_{j_0}-x_{j_0+1}| < 1/N$, or else 
\item for at least one $1 \leq j_0 \leq k-1$ we have $x_{j_0} \in [y_{m_{j_0}},y_{m_{j_0}+1}] \bk I_{m_{j_0}}$, i.e., $x_{j_0}$ lies in an interval of length $(ZN)^{-2k}$.
\end{enumerate}
We thus see that
\begin{align*}
&\int_{\mc{B}(Z) \bk \mc{B}_{N}(Z)} e^{-\frac{1}{2}\mbf{x}\cdot \mbf{x}} d\mbf{x} \\
&\leq \sum_{1 \leq i \leq k-1} \int_{\mc{B}(Z)} 1_{|x_i-x_{i+1}| < 1/N}(\mbf{x}) e^{-\frac{1}{2}\mbf{x} \cdot \mbf{x}} d\mbf{x} \\
&+ \sum_{1 \leq i \leq k-1} \sum_{\substack{\mbf{m} \in \{0,\ldots,2ZN-1\}^k\\ m_j < m_{j+1} \\ \forall  1\leq j \leq k-1}} \int_{B_{N,Z}(\mbf{m})} 1_{|x_i-y_{m_i+1}| < (NZ)^{-2k}}(\mbf{x}) e^{-\frac{1}{2}\mbf{x}\cdot \mbf{x}} d\mbf{x} \\
&\ll_k 1/N + (2ZN)^k \cdot (NZ)^{-2k} \ll_k 1/N.
\end{align*}
Since the boxes $\{\mc{B}_{N,Z}(\mbf{m}) : \mbf{m} \in \{0,\ldots,2NZ-1\}^k, m_j < m_{j+1} \forall 1 \leq j \leq k-1\}$ are mutually disjoint, the above estimates show that
$$
\sum_{\substack{\mbf{m} \in \{0,\ldots,2N Z-1\}^{k} \\ m_j < m_{j+1} \\ \forall 1 \leq j \leq k-1}} \int_{B_{N,Z}(\mbf{m})} e^{-\frac{1}{2}\mbf{x} \cdot \mbf{x}} d\mbf{x} = \int_{\mc{B}} e^{-\frac{1}{2}\mbf{x} \cdot \mbf{x}} d\mbf{x} + O_k\left(N^{-1} + e^{-Z^2/4}\right).
$$
and so
$$
(2\pi)^{-k/2}\sum_{\substack{\mbf{m} \in \{0,\ldots,2N Z-1\}^{k} \\ m_j < m_{j+1} \\ \forall 1 \leq j \leq k-1}} \int_{B_{N,Z}(\mbf{m})} e^{-\frac{1}{2}\mbf{x} \cdot \mbf{x}} d\mbf{x} = \frac{1}{k!} + O_k\left(1/N + e^{-Z^2/4}\right).
$$
We now suppose that $x$ is sufficiently large relative to $Z$ and $N$ (which, themselves, are chosen large enough in terms of $k$). By Proposition \ref{prop:EKCusp}, it follows that
\begin{align}
&\sum_{\substack{\mbf{m} \in \{0,\ldots,2N Z-1\}^{k} \\ m_j < m_{j+1} \\ \forall 1 \leq j \leq k-1}} X^{-1}|\{n \leq x : n+a_j \in \mc{N}, X_j(n+a_j) \in I_{m_j} \text{ for all } 1 \leq j \leq k\}| \label{eq:truncBoxes}\\
&= (2\pi)^{-k/2}\sum_{\substack{\mbf{m} \in \{0,\ldots,2N Z-1\}^{k} \\ m_j < m_{j+1} \\ \forall 1 \leq j \leq k-1}}  \int_{B_{N,Z}(\mbf{m})} e^{-\frac{1}{2}\mbf{x} \cdot \mbf{x}} d\mbf{x} + o_k(1) \nonumber\\
&= \frac{1}{k!} + O_k\left(1/N + e^{-Z^2/4}\right). \nonumber
\end{align}
Finally, we relate this to $T_{\mbf{a}}(x)$. Note that if $\sqrt{x} < n \leq x$ is counted by the sum in \eqref{eq:truncBoxes} then there is some $\mbf{m} \in \{0,\ldots, 2NZ-1\}^k$ with $m_j < m_{j+1}$ such that for all $1 \leq j \leq k-1$,
\begin{align*}
X_j(n+a_j) &\leq y_{m_j+1} - (NZ)^{-2k} \leq y_{m_{j+1}} - (NZ)^{-2k} \\
&\leq X_j(n+a_{j+1})-1/(NZ)^{2k} < X_j(n+a_{j+1}) - \frac{1}{(\log\log x)^{1/3}},
\end{align*}
provided that $x$ is sufficiently large in terms of $N$ and $Z$. We find, therefore, that
\begin{align*}
X^{-1}T_{\mbf{a}}(x) &\geq \sum_{\substack{\mbf{m} \in \{0,\ldots,2N Z-1\}^{k} \\ m_j < m_{j+1} \\ \forall 1 \leq j \leq k-1}} X^{-1}|\{n \leq x : n+a_j \in \mc{N}, X_j(n+a_j) \in I_{m_j} \text{ for all } 1 \leq j \leq k\}| \\
&= \frac{1}{k!} + O_k(1/N+e^{-Z^2/4}),
\end{align*}
for large enough $x$, and hence as $x \ra \infty$ we obtain
$$
X^{-1}T_{\mbf{a}}(x) \geq \frac{1}{k!} + O_k(1/N+e^{-Z^2/4}).
$$
For each $\sg \in S_k$, write $\mbf{a}_{\sg} = (a_{\sg(1)},\ldots,a_{\sg(k)})$. By the above argument, applied to each $\mbf{a}_{\sg}$ (there are only finitely many of these, so we choose $x$ large enough according to each of these cases), it thus follows that
\begin{align*}
X^{-1}T_{\mbf{a}_{\sg}}(x) \geq \frac{1}{k!} + O_k\left(1/N + e^{-Z^2/4}\right) \text{ for all } \sg \in S_k.
\end{align*}
We now claim that in fact $X^{-1}T_{\mbf{a}}(x) = 1/k! + o(1)$. Indeed, suppose for the sake of contradiction that there is an infinite sequence $\{x_l\}_l$ such that
$$
X(x_l)^{-1}T_{\mbf{a}_{\sg}}(x_l) \geq 1/k! + \e,
$$
for some $\e > 0$, as $l \ra \infty$. Then by the above, we obtain, for $x$ large enough,
\begin{align*}
\sum_{\sg \in S_k} T_{\mbf{a}_{\sg}}(x_l) &\geq \frac{X(x_l)}{k!} +X(x_l) \e/2 + X(x_l)(k!-1) \cdot\left(1/k! + O_k(1/N + e^{-Z^2/4})\right) \\
&= X(x_l)(1+\e/2 + O_k(1/N+e^{-Z^2/4})).
\end{align*} 
On the other hand, we clearly have
$$
\sum_{\sg \in S_k} |\{n\leq x : n+a_{\sg(j)} \in \mc{N}, X_j(n+a_{\sg(j)}) < X_{j+1}(n+a_{\sg(j+1)})  \ \forall 1 \leq j \leq k-1\}| \leq X(x_l),
$$
This produces a contradiction whenever $Z,N$ are large enough relative to $\e$, which is permissible as long as $x$ is chosen sufficiently large in terms of these. We thus conclude that for any $\e > 0$ we can find $x$ large enough (after picking $N$ and $Z$ suitably large in terms of $\e$) such that
$$
\frac{1}{k!} - \e \leq X^{-1}|\{n \leq x: n+a_j \in \mc{N}, |b_1(n+a_{1})| < \cdots < |b_k(n+a_k)|\}| \leq \frac{1}{k!} + \e,
$$
and the proof is complete.
\end{proof}
\begin{proof}[Proof of Theorem \ref{thm:condSigns}]
Let us say that a sign pattern $\mbf{\eps} \in \{0,1\}^k$ is \emph{good} if there is an integer $0 \leq r \leq k$ such that 
$$
\eps_i = \begin{cases} -1 &\text{ if $1 \leq i \leq r$,} \\ +1 &\text{ if $r+1 \leq i \leq k$} \end{cases}
$$
(with the obvious convention if $r = 0$ or $r = k$). Suppose $n$ satisfies $n+a_j \in \mc{N}$, and define $\sg_j(n+a_j) := \text{sign}(b_j(n+a_j))$ for all $1 \leq j \leq k$. Observe that
$$
b_1(n+a_j) < \cdots < b_k(n+a_k) \Rightarrow (\sg_1(n+a_1),\ldots,\sg_k(n+a_k)) \text{ is a good sign pattern:}
$$
indeed, if not then there is a $1 \leq j \leq k-1$ such that $\sg_j(n+a_j) = +1$ and $\sg_{j+1}(n+a_{j+1}) = -1$, but then by assumption we have $0 < |b_j(n+a_j)| < -|b_{j+1}(n+a_{j+1})|$, an obvious contradiction. Letting $\mc{G} \subseteq \{-1,+1\}^k$ denote the set of good sign patterns, we thus obtain 
\begin{align*}
&|\{n \leq x : b_1(n+a_1) < \cdots < b_k(n+a_k)\}| \\
&\geq |\{n \leq x : n+a_1,\ldots,n+a_k \in \mc{N}, b_1(n+a_j) < \cdots < b_k(n+a_k)\}| \\
\begin{split}
&= \sum_{\mbf{\eps} \in \mc{G}} |\{n \leq x : n+a_j \in \mc{N}, \sg_j(n+a_j) = \e_j \forall 1 \leq j \leq k \\
&\text{ and } \eps_i |b_i(n+a_i)| < \eps_j|b_j(n+a_j)| \text{ for all $1 \leq i < j \leq k$}\}| 
\end{split}
\\
&=: \sum_{\mbf{\eps} \in \mc{G}} T_{\mbf{a},\mbf{\eps}}(x).
\end{align*}
We observe that if $n$ is counted by $T_{\mbf{a},\mbf{\eps}}(x)$ then the inequality $b_i(n+a_i) < b_j(n+a_j)$ with $1 \leq i < j \leq k$ is vacuous for $\mbf{\eps} \in \mc{G}$ whenever $\e_i \e_j = -1$ (as in this case $\eps_i = -1 = -\eps_j$ and the inequality $-|b_i(n+a_i)| < |b_j(n+a_j)|$ is trivially always valid). On the other hand, if $1 \leq i < j \leq k$ then we also note that
\begin{align*}
\eps_i = \eps_j = +1 &\Leftrightarrow |b_i(n+a_i)| < |b_j(n+a_j)| \\
\eps_i = \eps_j = -1 &\Leftrightarrow |b_i(n+a_i)| > |b_j(n+a_j)|. 
\end{align*}
What transpires is that if $\eps \in \mc{G}$ then $n$ is counted by $T_{\mbf{a},\mbf{\eps}}(x)$ if, and only if, $(|b_1(n+a_1)|,\ldots,|b_k(n+a_k)|)$ is contained in the set
$$
\mc{A}_{\eps} := \{\mbf{x} \in \mb{R}^k : \eps_i = \eps_j = +1 \Rightarrow x_i < x_j, \eps_i = \eps_j = -1 \Rightarrow x_i > x_j \text{ for all $1 \leq i < j \leq k$}\}.
$$
Given $\mbf{\eps} \in \mc{G}$, let us denote by $p(\mbf{\eps})$ the integer $0 \leq r \leq k$ such that $\eps_i = -1$ for all $1 \leq i \leq r$ (with $p((+1,\ldots,+1)) := 0$); trivially, for each $0 \leq r\leq k$ there is a unique good sign pattern $\mbf{\eps}$ with $p(\mbf{\eps}) = r$.  We observe directly that if $p(\mbf{\eps}) = r$ then we can decouple the first $r$ variables $x_1,\ldots,x_r$ from the last $k-r$ variables in the Gaussian integral restricted to $\mc{A}_{\mbf{\eps}}$, giving
\begin{align*}
(2\pi)^{-k/2} \int_{\mc{A}_{\mbf{\eps}}} e^{-\frac{1}{2}\mbf{x} \cdot \mbf{x}} d\mbf{x} &= (2\pi)^{-r/2} \left(\int_{\mb{R}^r} e^{-\frac{1}{2}(x_1^2+\cdots + x_r^2)} 1_{x_1 > \cdots > x_r} dx_1\cdots dx_r\right)\\
&\quad \cdot (2\pi)^{-(k-r)/2} \left(\int_{\mb{R}^{k-r}} e^{-\frac{1}{2}(x_{r+1}^2+\cdots + x_k^2)} 1_{x_{r+1} < \cdots < x_k} dx_{r+1}\cdots dx_k\right) \\
&= \frac{1}{r!(k-r)!}
\end{align*}
via a symmetry argument as in the deduction of Theorem \ref{thm:genCusp}. Let us therefore note the identity
\begin{align} \label{eq:intIden}
2^{-k} \sum_{\mbf{\eps} \in \mc{G}} (2\pi)^{-k/2} \int_{\mc{A}_{\mbf{\eps}}} e^{-\frac{1}{2}\mbf{x} \cdot \mbf{x}} d\mbf{x} = \frac{1}{k!} \cdot 2^{-k} \sum_{0 \leq r \leq k} \binom{k}{r} = \frac{1}{k!},
\end{align}
which we will use shortly. \\
We use a covering argument as in the deduction of Theorem \ref{thm:genCusp}; since the ideas are similar we merely sketch the argument in this case. Letting $N,Z \geq 1$ be large parameters and $\mbf{\eps} \in \mc{G}$, we define 
\begin{align*}
\mc{M}_{\mbf{\eps}} := \{\mbf{m} \in \{0,\ldots,2NZ-1\}^k : \ &\eps_i = \eps_j = +1 \Rightarrow m_i < m_j, \\
&\eps_i=\eps_j = -1 \Rightarrow m_i > m_j \text{ for all $1 \leq i < j \leq k$}\},
\end{align*}
and set 
$$
B_{N,Z}(\mbf{m}) = \prod_{1 \leq j \leq k} I_{m_j} = \prod_{1 \leq j \leq k} [y_{m_j},y_{m_j+1}-(NZ)^{-2k}]
$$
with $y_l := -Z + l/N$, as before. Arguing as in the previous deduction, we end up with
$$
\sum_{\mbf{m} \in \mc{M}_{\mbf{\eps}}} (2\pi)^{-k/2} \int_{B_{N,Z}(\mbf{m})} e^{-\frac{1}{2}\mbf{x} \cdot \mbf{x}} d\mbf{x} = (2\pi)^{-k/2} \int_{\mc{A}_{\mbf{\eps}}} e^{-\frac{1}{2}\mbf{x} \cdot \mbf{x}} d\mbf{x} + O_k(1/N + e^{-Z^2/4}).
$$
We assume Conjecture \ref{conj:EllIntro}. By Theorem \ref{thm:condSigns}, for each $\mbf{m} \in \mc{M}_{\mbf{\eps}}$ we obtain
\begin{align*}
&X^{-1}|\{n \leq x : n + a_j \in \mc{N}, \sg_j(n+a_j) = \eps_j, X_j(n+a_j) \in I_{m_j} \text{ for all } 1 \leq j \leq k\}| \\
&= 2^{-k}(2\pi)^{-k/2} \int_{B_{N,Z}(\mbf{m})} e^{-\frac{1}{2}\mbf{x} \cdot \mbf{x}} d\mbf{x} + o_k(1),
\end{align*}
so taking $x$ large enough relative to $N$ and $Z$, we obtain
\begin{align*}
&\sum_{\mbf{\eps} \in \mc{G}} \sum_{\mbf{m} \in \mc{M}_{\mbf{\eps}}} X^{-1} |\{n \leq x : n+a_j, \sg_j(n+a_j) = \eps_j, X_j(n+a_j) \in I_{m_j} \text{ for all } 1 \leq j \leq k\}| \\
&= 2^{-k} \sum_{\mbf{\eps} \in \mc{G}} (2\pi)^{-k/2} \int_{\mc{A}_{\mbf{\eps}}} e^{-\frac{1}{2}\mbf{x} \cdot \mbf{x}} d\mbf{x} + O_k(1/N+e^{-Z^2/4}) = \frac{1}{k!} + O_k(1/N+e^{-Z^2/4}),
\end{align*}
Finally, arguing as in the previous deduction, we may obtain
\begin{align*}
\begin{split}
T_{\mbf{a},\mbf{\eps}}(x) \geq &|\{\sqrt{x} < n \leq x : n+a_j \in \mc{N}, \sg_j(n+a_j) = \eps_j \text{ for all $1 \leq j \leq k$,} \\
 &\eps_i = \eps_j = +1 \Rightarrow X_i(n+a_i)< X_j(n+a_j) - (\log\log x)^{-1/3}, \\
 &\eps_i = \eps_j = -1 \Rightarrow X_i(n+a_j) > X_j(n+a_j) + (\log\log x)^{-1/3}) \text{ for all } 1 \leq i < j \leq k\}| \\
&+ O(\sqrt{x})
\end{split}
\\
&\geq \sum_{\mbf{m} \in \mc{M}_{\mbf{\eps}}} |\{n \leq x : n+a_j, \sg_j(n+a_j) = \eps_j, X_j(n+a_j) \in I_{m_j} \text{ for all } 1 \leq j \leq k\}| + O(\sqrt{x}),
\end{align*}
for $x$ large enough, so combined with the previous estimates, it follows that
$$
\sum_{\mbf{\eps} \in \mc{G}} T_{\mbf{a},\mbf{\eps}}(x) \geq X\left(\frac{1}{k!} + O_k(1/N + e^{-Z^2/4})\right),
$$
if $x$ is large enough. Finally, by the permutation trick employed at the end of the deduction of Theorem \ref{thm:genCusp}, we may infer immediately that as $x \ra \infty$,
$$
\sum_{\mbf{\eps} \in \mc{G}} T_{\mbf{a},\mbf{\eps}}(x) = X\left(\frac{1}{k!} + o_k(1)\right),
$$
whence we deduce that
\begin{align*}
&|\{n \leq x : b_1(n+a_1) < \cdots < b_k(n+a_k)\}| \\
&\geq |\{n \leq x : n+a_1,\ldots,n+a_k \in \mc{N}, b_1(n+a_1) < \cdots < b_k(n+a_k)\}| \\
&= X\left(\frac{1}{k!} + o_k(1)\right),
\end{align*}
as claimed. 
\end{proof}

\begin{proof}[Proof of Theorem \ref{thm:uncondSigns}]
It follows from partial summation that for any given set $S \subset \mb{N}$ one has
$$
\bar{d}(S) \geq \limsup_{Y \ra \infty} \frac{1}{\log Y} \sum_{n \leq Y} \frac{1_S(n)}{n}.
$$
It thus suffices to establish that the limit on the RHS of this expression is non-zero. Set now
\begin{align*}
S_{\mbf{a}} &:= \{n \in \mb{N} : b_1(n+a_1) < \cdots < b_k(n+a_k)\} \\
S_{\mbf{a},\mbf{\eps}} &:= \{n \in \mb{N}: n+a_j \in \mc{N}, \sg_j(n+a_j) = \eps_j, \eps_i |b_i(n+a_i)| < \eps_j|b_j(n+a_j)| \forall 1 \leq i < j \leq k\},
\end{align*}
with either $k = 2$, or $k = 3$ and at least two of $b_1,b_2, b_3$ the same. Arguing as in the deduction of Theorem \ref{thm:condSigns} (with logarithmic averages instead of C\'{e}saro ones, by partial summation), we obtain from Theorem \ref{thm:uncondSigns} that
\begin{align*}
\frac{1}{\log Y} \sum_{n \leq Y} \frac{1_{\mc{S}_{\mbf{a}}}(n)}{n} &\geq \sum_{\mbf{\eps} \in \mc{G}} \frac{1}{\log Y}\sum_{\substack{n \leq Y \\ n+a_j \in \mc{N} \\ \forall 1 \leq j \leq k}} \frac{1_{\mc{S}_{\mbf{a},\mbf{\eps}}}(n)}{n} = \left(\frac{1}{k!}+o_k(1)\right) \frac{1}{\log Y}\sum_{\substack{n \leq Y \\ n+a_j \in \mc{N} \\ \forall 1 \leq j \leq k}} \frac{1}{n},
\end{align*}
Taking limits, we obtain the desired conclusions since the (normalized) logarithmic sum here is $\gg_{k, \mbf{a}, \mbf{I}} 1$ as $Y \ra \infty$. 
\end{proof}

\bibliography{tau1}
\bibliographystyle{plain}
\nocite{Kub}
\end{document}